\documentclass[sn-basic,Numbered]{sn-jnl}
\usepackage{geometry}
\usepackage{multirow}%
\usepackage{amsmath,amssymb,amsfonts}%
\usepackage{amsthm}%
\usepackage{mathrsfs}%
\usepackage[title]{appendix}%
\usepackage{xcolor}%
\usepackage{textcomp}%
\usepackage{manyfoot}%
\usepackage{booktabs}%
\usepackage{algorithm}%
\usepackage{algorithmicx}%
\usepackage{algpseudocode}%
\usepackage{listings}%
\usepackage{amssymb}
\usepackage{url}
\usepackage{bm}
\usepackage{xspace}
\newtheorem{theorem}{Theorem}
\newtheorem{remark}{Remark}%

\DeclareMathAlphabet{\mathitbf}{OML}{cmm}{b}{it}
\def\jmath{j}

\def\jmath{j}

\usepackage[labelformat=simple]{subcaption}

\newcommand{\frwidth}{\epsilon}

\newcommand{\EXP}[1]{\mathbb{E}\left[#1\right ]}

\newcommand{\Var}[1]{\mathbb{V}\left[#1\right ]}

\def\min{\mbox{min}}

\def\sol{c}

\def\D{\mathcal{D}}

\newcommand{\MSE}{\text{MSE}}

\def\xib{\bm{\xi}}

\def\bx{\mathbf{x}}

\newcommand{\conc}{c} 
\newcommand{\pres}{p} 
\newcommand{\poro}{\phi} 
\newcommand{\perm}{{\mathbf{K}}} 
\newcommand{\dens}{\rho} 
\newcommand{\visc}{\mu} 
\newcommand{\dvel}{{\mathbf{q}}} 
\newcommand{\grav}{{\mathbf{g}}} 
\newcommand{\disp}{{\mathbf{D}}} 

\newtheorem{defn}[theorem]{Definition}

\begin{document}

\title[]{ Estimation of uncertainties in the density driven flow in fractured porous media using MLMC}



\author[1]{\fnm{Dmitry} \sur{Logashenko}}\email{dmitry.logashenko@kaust.edu.sa}

\author*[2]{\fnm{Alexander} \sur{Litvinenko}}\email{litvinenko@uq.rwth-aachen.de}

\author[1,2]{\fnm{Raul} \sur{Tempone}}\email{raul.tempone@kaust.edu.sa}

\author[1]{\fnm{Gabriel} \sur{Wittum}}\email{gabriel.wittum@kaust.edu.sa}

\affil[1]{\orgdiv{CSE Department}, \orgname{KAUST}, \city{Thuwal-Jeddah}, \country{Saudi Arabia}}

\affil[2]{\orgdiv{Department of Mathematics}, \orgname{RWTH Aachen}, \orgaddress{\city{Aachen}, \country{Germany}}}


\abstract{
We use the Multi Level Monte Carlo method to estimate uncertainties in a Henry-like salt water intrusion problem with a fracture. The flow is induced by the variation of the density of the fluid phase, which depends on the mass fraction of salt. While the fracture's location is fixed, its aperture is uncertain. 
In our setting, porosity and permeability vary spatially and recharge is time-dependent. 
So we introduce three random variables, one controlling both the porosity and permeability fields, one for the fracture width and one for the intensity of recharge.
For each realization of these uncertain parameters, the evolution of mass fraction and pressure fields is modeled using a system of non-linear, time-dependent PDEs with a solution discontinuity at the fracture. These uncertainties propagate, affecting the distribution of salt concentration, a key factor in water resource quality. We show that the MLMC method can be successfully applied to this problem. It significantly reduces the computational cost compared to classical Monte Carlo methods by effectively balancing discretisation and statistical errors, and by evaluating multiple scenarios over different spatial and temporal mesh levels. The deterministic PDE solver, using the ug4 library, runs in parallel to compute all stochastic scenarios.}



\maketitle
\tableofcontents
\section{Introduction}\label{sec1}

Density-driven flow is a phenomenon observed in several natural and engineering contexts, including weather patterns, atmospheric and oceanic circulation, geothermal energy extraction, and cooling of electronic and nuclear reactors. In hydrogeology, modelling density-driven subsurface flow and salt transport in fractured aquifers presents significant challenges. These arise from the complex, parameter-sensitive flow patterns caused by fractures and various uncertain parameters. Fractures introduce heterogeneity into the aquifer, affecting fluid flow and making it difficult to predict key properties. In addition, fractures can become sealed over time, reducing their permeability and further affecting fluid flow. 

Uncertainties in fracture geometry, porosity, permeability and recharge have a significant impact on the prediction of the evolution of subsurface salinity (cf. \cite{Angot, Bastian2000, Reichenberger+Helmig, Martinez2006, Sorek2001}). In real reservoirs, the number, location and aperture of fractures are inherently uncertain. One possible solution is to model these uncertainties, build a surrogate model and perform sensitivity analysis (see \cite{Sudret_sparsePCE, sudret2008global}) to identify critical uncertainties and assess their impact on the solution. However, this approach is complex. 

The combination of MLMC with the geometric multigrid method is a promising strategy. The geometric multigrid method (GMG) uses a hierarchy of grids and Multi-Level Monte Carlo (MLMC) involves sampling from different grids. The multigrid method is already implemented in our UG4 software library and in many other codes.
We note that both MLMC and surrogate models can be implemented in an almost black-box fashion. 

In our previous work \cite{LOGASHENKO2024_JCP}, we applied the MLMC method to a similar problem, albeit without a fracture in the aquifer. The problem addressed in this paper is more complex, involving a greater number of equations, additional uncertainties, and a discontinuous solution. In other related works \cite{LitLog3D_20} and \cite{Litvinenko-UQ-2021}, the authors employed the Polynomial Chaos Expansion (PCE) method to solve the Elder problem, another density-driven flow problem. However, the Elder problem differs significantly from the current problem. While the authors successfully applied PCE, they encountered certain limitations, and it remains unclear whether their methods and insights can be effectively extended to the present problem with a fractured reservoir. Multi-level multi-fidelity surrogate models have been explored in \cite{Multifidelity19}.

In \cite{Riva2015}, the authors are solving the seawater intrusion problem under multiple sources of uncertainty. They use the generalised polynomial chaos expansion approximation (gPCE) to investigate how incomplete knowledge of system properties affects the assessment of global quantities. They consider very different input uncertainties, such as the gravity number, the permeability anisotropy ratio, and the transverse and longitudinal Peclet numbers. Their uncertainties are scalar values, whereas our uncertainties include random fields (porosity and permeability) and random process (recharge).

In \cite{GRILLO2010,GLLSW-JPM2012,Reiter14}, the authors modeled and computed the density-driven flow in fractured porous media with deterministic settings. The fractures considered were represented by a $(d - 1)$-dimensional manifolds in a $d$-dimensional domain. They developed a special numerical technique based on the finite volume method and implemented it in the ug3 library. In our present work, we apply the same model, the same discretization and solvers to compute the multiple uncertain scenarios. Other methods based on $(d-1)$-dimensional representation of fractures were considered in \cite{Angot,Bastian2000,Reichenberger+Helmig,Martinez2006,Sorek2001,Shikaze1998,GrafTherrien2007}. It is worth to remark that this type of approaches with the low-dimensional representation of the fractures becomes inaccurate for large fracture apertures. As demonstrated in \cite{GLLSW-JPM2012} where simulation results of similar deterministic settings with $d$- and $(d-1)$-dimensional representations of the fractures are compared, this issue becomes significant for the density-driven flow: Due to the non-linearity of the equations, vortices arising in thick fractures change the flow and transport in the entire domain essentially. Taking this observation into account, in the present work, we restrict ourselves to sufficiently thin fractures.

Many techniques can be used to quantify uncertainties. One classical method is Monte Carlo (MC) sampling. It has a well-known disadvantage --- a slow convergence $\mathcal{O}(\frac{1}{\sqrt{N}})$. Other relative recent techniques such as surrogate models and stochastic collocation may require a few hundred time-consuming simulations and assume a certain smoothness of the quantity of interest (QoI).
Another class of methods is the class of perturbation methods \cite{CREMER15_Fingers}. The idea is to decompose the QoI with respect to random parameters in a Taylor series. The higher order terms can be neglected for small perturbations, simplifying the analysis and numerics. These methods assume that the random perturbations are small. For larger perturbations, these methods usually do not work. 


Recent efforts have benchmarked various discrete fracture models, as detailed in \cite{Flemisch_IWR_2018} and \cite{Berre_2021}, albeit in deterministic settings rather than uncertain ones. In \cite{Flemisch_IWR_2018}, several benchmarks were presented for numerical schemes addressing single-phase fluid flow in fractured porous media. These benchmarks compared multiple solution strategies and progressively increased the complexity of network geometry—such as intersecting fractures—and physical parameters, including low and high fracture-matrix permeability ratios and heterogeneous fracture permeabilities. Similarly, \cite{Berre_2021} introduced four benchmark cases for single-phase flow in three-dimensional fractured porous media. These cases were meticulously designed to test the capability of various methods to manage the complexities inherent in the geometric structures of fracture networks. A total of 17 numerical methods were collected and compared.

As a model problem, we consider the benchmark from \cite{GRILLO2010} which is a generalization of the Henry problem first introduced in \cite{henry1964effects}. The Henry problem became a benchmark for numerical solvers for the density-driven groundwater flow (see \cite{Voss_Souza,Simpson04_Henry,Simpson2003,Dhal_review22}.

\textbf{Novelty.} We combined the geometric multi-grid and MLMC methods to quantify uncertainties in a discontinuous, time-dependent, nonlinear, density-driven flow problem characterized by a system of coupled partial differential equations. The primary objective was to apply and to test the MLMC method. Both the MLMC method and parallel geometric multigrid methods were employed to estimate the mean and variance of the salt mass fraction.
The computational domain encompasses a fracture. Three parameters —-- porosity, recharge and fracture aperture --- are treated as uncertain. The porosity depends on the spatial coordinates, the recharge is time-dependent, and the fracture aperture is modeled by a random variable. The permeability of the medium is a function of its porosity. Our numerical experiments confirm that, under certain settings, the MLMC method is applicable and can compute the mean up to 100 times faster than standard Monte Carlo methods.

This work is structured as follows. Section~\ref{sec:Model} describes the model. Section~\ref{sec:NumMethods} provides a description of the numerical methods, with a review of the well-known MLMC method in Section~\ref{sec:MLMC}. Section~\ref{sec:numerics} 
presents the numerical results, including the analysis of the model problem, computation of various statistics, and the performance evaluation of both the MLMC method and the parallel multigrid solver. Finally, we conclude this work with a discussion in Section~\ref{sec:Conclusion}.\\

\section{Modeling and Problem Settings}
\label{sec:Model}

We examine the density-driven flow of a liquid phase (salt solution in water) within a fractured aquifer, where the salt mass fraction varies. Fractures in the aquifer are geological formations with a negligible size in one geometric direction compared to the domain's scales, but the permeability of their filling material is significantly higher than that of the surrounding medium. This disparity presents substantial challenges when treating fractures as full-dimensional subdomains in simulations. Specifically, resolving these features requires very fine grids, and the sharp contrast in model parameters drastically reduces numerical method efficiency. Consequently, several approaches have been developed to address these issues (cf. \cite{Angot, Reichenberger+Helmig, Fumagalli+Scotti-2013, dAngelo+Scotti, Martin+Jaffre, Hoteit2008891, Reiter14}).

A notable advantage of these models for sampling methods is the ability to alter the fracture width without remeshing the entire domain. In this study, we assume the porous matrix remains immobile, and the fractures retain their positions.
\subsection{Governing equations for the flow}
\label{subsec:PDEs}

In this section, we briefly describe the partial differential equations (PDEs) of the model. For further details, we refer to \cite{Reiter2012,Reiter14}. We denote by $\D \subseteq \mathbb{R}^d$ the entire domain representing the aquifer with $\mathscr{M} \subseteq \mathscr{D}$ being its part surrounding the fracture. The porosity of this part is $\phi_m : \mathscr{M} \to (0,1)$ and its permeability $\perm_m : \mathscr{M} \to \mathbb{R}^{d \times d}$.

The fluid phase in $\mathscr{M}$ is characterized by the salt mass fraction $\conc_m (t, \mathbf{x}): [0, +\infty) \times \mathscr{M} \to [0, 1]$ (where $1$ corresponds to the saturated brine) and pressure $\pres_m (t, \mathbf{x}): [0, +\infty) \times \mathscr{M} \to \mathbb{R}$. Thus, the flow and transport in $\mathscr{M}$ obey the standard conservation of mass laws for the entire liquid phase and the salt
\begin{equation} \label {e_cont_eq_m}
 \left .
 \begin{array}{l}
  \partial_t (\poro_m \dens_m) +  \nabla \cdot (\dens_m \dvel_m) = 0 \\
  \partial_t (\poro_m \dens_m \conc_m) + \nabla \cdot (\dens_m \conc_m \dvel_m - \dens_m \disp_m \nabla \conc_m) = 0
 \end{array}
 \right \} \quad x \in \mathscr{M}.
\end{equation}
In \eqref{e_cont_eq_m}, we assume the Darcy's law for the velocity:
\begin{eqnarray} \label {e_Darcy_vel_m}
 \dvel_m = - \frac{\perm_m}{\visc_m} (\nabla \pres_m - \dens_m \grav), \qquad x \in \mathscr{M}.
\end{eqnarray}
In (\ref{e_cont_eq_m}--\ref{e_Darcy_vel_m}), $\dens_m = \dens (\conc_m)$ and $\visc_m$ are the density and the viscosity of the liquid phase, $\disp_m (t, \mathbf{x}, \dvel_m)$ denotes the diffusion and mechanical dispersion tensor.

The fractures are also considered to be filled with a porous medium. In particular, we assume Darcy's law inside the fractures. Without loss of generality, we consider only one fracture in this work, see \cite{Reiter2012,Reiter14} for a generalization. 
The fracture is represented by a surface $\mathscr{S} \subset \D$, $\overline{\mathscr{M}} \cup \overline{\mathscr{S}} = \overline{\D}$, $\mathscr{M} \cap \mathscr{S} = \emptyset$. However, we distinguish between the two sides of the fracture, $\mathscr{S}^{(1)}$ and $\mathscr{S}^{(2)}$, which geometrically coincide with $\mathscr{S}$, but virtually are interfaces of the fracture with $\mathscr{M}$ and have opposite normals $\mathbf{n}^{(1)}$ and $\mathbf{n}^{(2)}$.

For the fluid phase inside the fracture, let $\conc_f: [0, +\infty) \times \mathscr{S} \to [0,1]$ and $\pres_f: [0, +\infty] \times \mathscr{S} \to \mathbb{R}$ are the salt mass fraction and the pressure averaged along the vertical, cf. \cite{GLLSW-JPM2012}. Note that in general the values of $\conc_f$ and $\pres_f$ are not equal to the limits of $\conc_m$ and $\pres_m$ at $\mathscr{S}$. Furthermore, the latter limits may be different on $\mathscr{S}^{(1)}$ and $\mathscr{S}^{(2)}$. We denote them by $\conc_m^{(1)}$ and $\conc_m^{(2)}$, as well as $\pres_m^{(1)}$ and $\pres_m^{(2)}$ as functions on $\mathscr{S}^{(1)}$ and $\mathscr{S}^{(2)}$, respectively. Two sets of equations are given for the fracture: the laws describing the flow along the fracture and those modeling the mass exchange through the interfaces $\mathscr{S}^{(1)}$ and $\mathscr{S}^{(2)}$.

For the flow and transport along the fracture, we impose the analogues of the mass conservation laws \eqref{e_cont_eq_m}:
\begin{equation} \label {e_cont_eq_f}
 \left .
 \begin{array}{l}
  \partial_t (\poro_f \frwidth \dens_f) +  \nabla^\mathscr{S} \cdot (\frwidth \dens_f \dvel_f) + Q_{fn}^{(1)} + Q_{fn}^{(2)} = 0 \\
  \partial_t (\poro_f \frwidth \dens_f \conc_f) + \nabla^\mathscr{S} \cdot (\frwidth \dens_f \conc_f \dvel_f - \frwidth \dens_f D_f \nabla^\mathscr{S} \conc_f) + P_{fn}^{(1)} + P_{fn}^{(2)} = 0
 \end{array}
 \right \} \quad x \in \mathscr{S},
\end{equation}
where $\nabla^\mathscr{S}$ denotes the differential operators on the manifold $\mathscr{S}$ and $\frwidth$ the fracture width. The tangential Darcy's velocity in the fracture is
\begin{eqnarray} \label {e_Darcy_vel_f}
 \dvel_f = - \frac{K_f}{\visc_f} (\nabla^\mathscr{S} \pres_f - \dens_f \grav),
  \qquad x \in \mathscr{S}.
\end{eqnarray}
In (\ref{e_cont_eq_f}--\ref{e_Darcy_vel_f}), $\rho_f = \rho (\conc_f)$ and $\mu_f = \mu (\conc_f)$ are the density and the viscosity of the fluid phase in the fracture. Furthermore, $\poro_f: \mathscr{S} \to (0, 1)$ denotes the porosity associated with the material in the fracture, $K_f : \mathscr{S} \to \mathbb{R}$ its tangential permeability and $D_f (t, x) : [0, +\infty) \times \mathscr{S} \to \mathbb{R}$ the tangential effective diffusion and dispersion tensor in the fracture.

The terms $Q_{fn}^{(k)}$ and $P_{fn}^{(k)}$, $k \in \{ 1,2 \}$, are the mass fluxes through the faces $\mathscr{S}^{(k)}$ of the fracture, i.e. the normal fluxes. They are defined as
\begin{equation} \label {e_cont_norm_flux}
 \left .
 \begin{array}{rl}
     Q_{fn}^{(k)} & := \rho (\conc_m^{(k)}) q_{fn}^{(k)}  \\
     P_{fn}^{(k)} & := \rho (\conc_m^{(k)}) c^{(k)}_\mathrm{upwind} q_{fn}^{(k)}
      - \rho (\conc_m^{(k)}) D_{fn}^{(k)} \dfrac{\conc_m^{(k)} - \conc_f}{\frwidth/2}
 \end{array}
 \right \} \quad x \in \mathscr{S}^{(k)},
\end{equation}
where $c^{(k)}_\mathrm{upwind} = \conc_m^{(k)}$ if $q_{fn}^{(k)} < 0$ but $c^{(k)}_\mathrm{upwind} = \conc_f$ if $q_{fn}^{(k)} \ge 0$ with
\begin{equation} \label {e_cont_norm_vel}
 q_{fn}^{(k)} := - \frac{K_{fn}^{(k)}}{\mu(\conc_f^{(k)})}
  \left [
   \frac{\pres_m^{(k)} - \pres_f}{\frwidth/2} - (\rho(\conc_f^{(k)}) - \rho_f) \grav \cdot \mathbf{n}^{(k)}
  \right ], \qquad x \in \mathscr{S}^{(k)}.
\end{equation}
In \eqref{e_cont_norm_vel}, $K_{fn}^{(k)} : \mathscr{S}^{(k)} \to \mathbb{R}$ is the normal permeability and $D_{fn}^{(k)} : [0, +\infty) \times \mathscr{S}^{(k)} \to \mathbb{R}$ is the normal effective diffusion and dispersion coefficient of the fracture -- bulk medium interface $\mathscr{S}^{(k)}$.

Systems \eqref{e_cont_eq_m} and \eqref{e_cont_eq_f} are coupled by the continuity of the mass fluxes of the interfaces $\mathscr{S}^{(k)}$, $k \in \{ 1, 2 \}$:
\begin{equation} \label {e_cont_eq_mf}
 \left .
 \begin{array}{rl}
   \rho (\conc_m^{(k)}) \dvel_m \cdot \mathbf{n}^{(k)} & = Q_{fn}^{(k)} \\
   \rho (\conc_m^{(k)}) \conc_m^{(k)} \dvel_m \cdot \mathbf{n}^{(k)}
    + \rho (\conc_m^{(k)}) \disp_m \nabla \conc_m \cdot \mathbf{n}^{(k)} & = P_{fn}^{(k)}
 \end{array}
 \right \}, \qquad x \in \mathscr{S}^{(k)}
\end{equation}
Furthermore, for the inner edges (in 3d) or the inner ends (in 2d) of the fractures, we impose the no-flux boundary conditions so that no mass exchange is possible there, cf. \cite{GRILLO2010,GLLSW-JPM2012}.

Equations \eqref{e_cont_eq_m}, \eqref{e_cont_eq_f} and \eqref{e_cont_eq_mf} with the velocities \eqref{e_Darcy_vel_m}, \eqref{e_Darcy_vel_f} and \eqref{e_cont_norm_vel} as well as particular specifications of the parameters $\rho$, $\mu$, $\poro_{m,f}$, $\perm_m$, $\disp_m$, $K_f$, $K_{fn}$, $D_f$ and $D_{fn}$ form the complete model for the density-driven flow in the porous medium with the fracture. This system must be closed by boundary conditions for $\conc_m$, $\conc_f$, $\pres_m$ and $\pres_f$ on $\partial \D$ as well as initial conditions for $\conc_m$ and $\conc_f$ at $t = 0$.

\subsection{Model problem settings}
\label{subsec_model_problem}

\begin{figure}[t!]
\begin{center}
\includegraphics [width=0.5\textwidth] {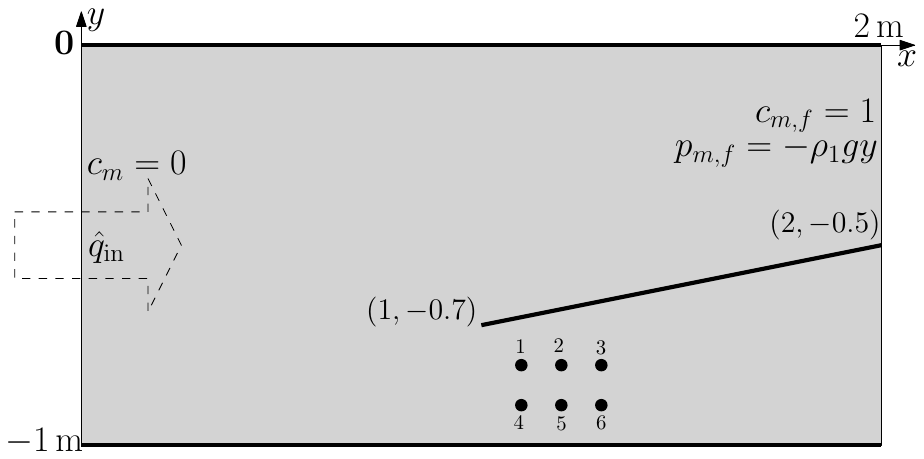}
\hfill
\includegraphics [width=0.45\textwidth] {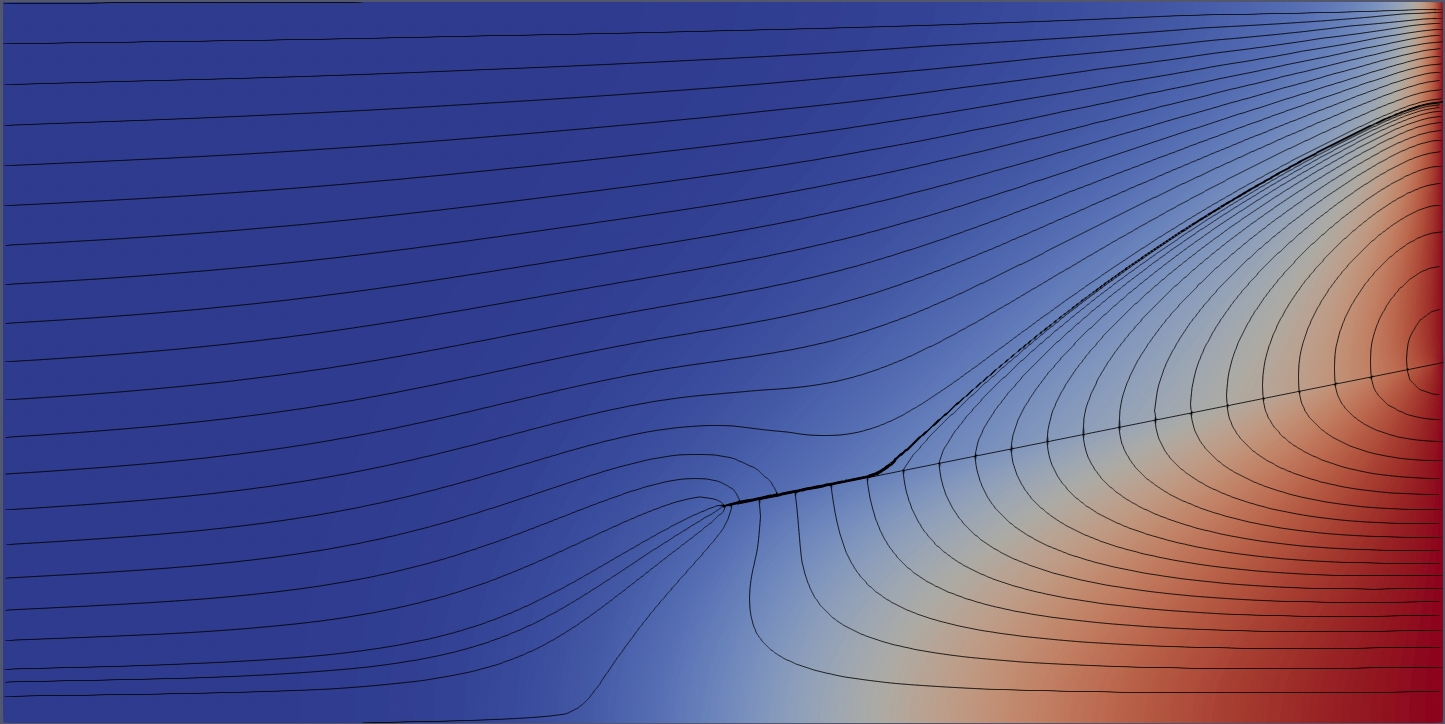}
\caption{Left: Computational geometry $[0,2]\times [-1,0]$, direction of recharge, boundary conditions, location of fracture (marked by black oblique line) and six points, solution in which we use as additional QoIs. Right: The flow streamlines (thin black lines) and the salt mass fraction. The dark red colour indicates a high salt mass fraction $\conc_m = 1.0$ or close to 1.0, and the dark blue corresponds to $\conc_m = 0$ (no salt). The thick black curved line indicates the area where the recharge flow meets the density difference induced flow.}
    \label{fig:Henry2d-scheme}
\end{center}    
\end{figure}

For our numerical tests, we choose a simple but very illustrative two-dimensional problem with one fracture, proposed in \cite{GRILLO2010}. It extends the setting of the Henry problem \cite{henry1964effects,Simpson04_Henry} with a fracture. The aquifer is represented by a rectangular domain $\D = [0, 2] \times [-1, 0]$ $[\mathrm{m}^2]$ completely saturated with the liquid phase, see Fig.~\ref{fig:Henry2d-scheme} (left). The fracture is located near the right ``sea side'' where the heavy salty water intrudes into the aquifer and cuts this boundary. The fracture coordinates are $(1,-0.7)$ and $(2,-0.5)$.

The system of the PDEs presented in Section \ref{subsec:PDEs} must be closed by specification of boundary conditions for $\conc_{m,f}$ and $\pres_{m,f}$, as well as initial conditions for $\conc_{m,f}$. In this work, we follow the settings from \cite{GRILLO2010} but add the uncertainty as described in Section \ref{subsec:PorosityVar} below. In particular, for the initial conditions, we set
\begin {equation}
 \left. \conc_m \right |_{t = 0} = 0, \qquad \left. \conc_f \right |_{t = 0} = 0.
\end {equation}

A scheme of the boundary conditions is shown in Fig.~\ref {fig:Henry2d-scheme} (left). On the right boundary, the Dirichlet conditions for $\conc_{m,f}$ and $\pres_{m,f}$ model the seawater intrusion:
\begin {equation} \label{e_seaside_bc}
	\left. \conc_{m,f} \right |_{x=2} = 1, \qquad \left . \pres_{m,f} \right |_{x=2} = - \rho_1 g y.
\end {equation}
On the left boundary, the inflow (recharge) of fresh water is imposed:
\begin {equation}
	\left. \conc_m \right |_{x=0} = 0, \qquad \left . \dens \dvel_m \cdot \mathbf{e}_x \right |_{x=0} = \hat{q}_{\mathrm{in}},
\end {equation}
where $\mathbf{e}_x = (1, 0)^\top$, and $\hat{q}_{\mathrm{in}}$ is a prescribed function of time, see Sect.~\ref{subsec:PorosityVar}. For the upper and the lower boundaries of $\D$, we impose no-flux boundary conditions for Equations \eqref{e_cont_eq_m}.

\begin{remark}
For the model problem presented, we rescale the mass fractions $\conc_{m,f}$ so that $\conc_{m,f} = 1$ corresponds to the most concentrated solution, i.e. seawater. Thus, seawater is considered as the abstract ``salt''. This explains the condition \eqref{e_seaside_bc} on the right boundary. Note that the density $\dens(c)$ is remapped accordingly, see \eqref{e_lin_density}. Therefore, the conservation equations for the salt (second equations in \eqref{e_cont_eq_m}, \eqref{e_cont_eq_f} and \eqref{e_cont_eq_mf}) are not affected.
\end{remark}

The flow and salt transport patterns that occur in this problem are discussed in \cite{GRILLO2010}. They are shown in Fig.~\ref{fig:Henry2d-scheme} (right). The colour codes the mass fraction $\conc_m$, with dark red corresponding to $\conc_m = 1$ and dark blue to $\conc_m = 0$, and generally $c_m\in[0,1]$. The lines are the streamlines of the flow. In the left part of $\D$, the flow is induced by the pure water injection through the left boundary. As in the Henry problem, the saltwater intruding through the lower part of the right boundary turns up and is flushed out through the upper part of the right side. In particular, the much more permeable fracture acts as a pathway for the strong washout that removes the saline water from the domain. For this purpose, the part of the domain below the fracture is mostly separated from the upper part of the domain --- as far as the capacity of the fracture allows. These phenomena depend on the width and permeability of the fracture.

In the case of time-independent boundary conditions, the mass fraction field converges to a steady state (cf. \cite{GRILLO2010}). However, prior to this there will be an initial period of instability which may last for a considerable time. This phase can even be seen for the time-dependent recharge $\hat{q}_{\mathrm{in}}$ considered below.

In this paper, as for the Henry problem (see \cite{Voss_Souza,Simpson04_Henry,Simpson2003}), we set
\begin {eqnarray} \label {e_lin_density}
 \dens (\conc) = \dens_0 + (\dens_1 - \dens_0) \conc, 
\end {eqnarray}
where $\dens_0$ is the density of the pure water and $\dens_1$ is the density of the brine (both considered as constants), and
$\visc = \text{const}$.
Besides that, we assume that porous medium is isotropic:
\begin {eqnarray} \label {e_scalar_perm}
 \perm_m = K_m \mathbf{I},
\end {eqnarray}
where $K_m : \mathscr{M} \to \mathbb{R}$. Furthermore, in this work, we neglect the mechanical dispersion, so that
\begin {eqnarray} \label {e_scalar_diff}
 \disp_m = \poro_m  D_0 \mathbf{I}, \qquad \disp_f = \poro_f D_0
\end {eqnarray}
where $D_0 \in \mathbb{R}$ is the constant effective diffusion coefficient of the salt in the liquid phase.
For the interface between the bulk medium and the fracture, we assume
\begin {eqnarray} \label {e_interface_values}
 K_{fn} = K_m, \qquad D_{fn} = \poro_m D_0, \qquad x \in \mathscr{S}.
\end {eqnarray}
Deterministic values of these parameters that are not varied in the scenarios are presented in Table~\ref {table-parameter}.

We examine the mass fraction $\conc_m$ at 6 locations where the largest variance is observed, i.e. we computed the variance in the whole domain and took locations where this variance is high. This corresponds to the situation when the water quality of wells reaching these points is controlled. The list of selected points follows (see Fig.~\ref{fig:Henry2d-scheme}):
\begin{align} \label{eq:6points} 
 \{\bx_i := (x_i,y_i)_{i=1,\ldots,6}\}
 = \{ &
	(1.1, -0.8),
    (1.2, -0.8),
    (1.3, -0.8), \nonumber \\
 &	(1.1, -0.9),
    (1.2, -0.9),
    (1.3, -0.9).
 \}
\end{align}
Whereas the other points can be considered, the reason to select these 6 points (and 6 small subdomains) is that not all other points are ``interesting'', i.e., not all points have significant variation in $\sol_m$.

\begin{table}[h!]
\footnotesize
\centerline
{
\begin{tabular}{|c|c|c|c|}
\hline
Symbol & Quantity & Value &Unit\\
\hline
	$D_\mathrm{mol}$ & Deterministic diffusion coefficient in the medium& $18.8571\cdot10^{-6}$ & [$\mathrm{m}^2 \textrm{ }\mathrm{s}^{-1}$] \\
	$\mathbf{g}$ & Gravity&$9.8$& [$\mathrm{m} \textrm{ }\mathrm{s}^{-2}$] \\
	$K_m$ & Deterministic permeability of the medium&$1.019368\cdot10^{-9}$& [$\mathrm{m}^2 $] \\
	$K_f$ & Deterministic permeability of the fracture&$1.019368\cdot10^{-6}$& [$\mathrm{m}^2 $] \\
		$\phi_m$ & Deterministic porosity of the medium&$0.35$& -\\
		$\phi_f$ & Deterministic porosity of the fracture&$0.7$& -\\
		$\mu$ & Viscosity&$10^{-3}$ & [$\mathrm{kg}\textrm{ } \mathrm{m}^{-1}\textrm{ } \mathrm{s}^{-1}$] \\
		$\rho_0$ & Density of water &$1\cdot10^{3}$& [$\mathrm{kg}\textrm{ } \mathrm{m}^{-3}$] \\
		$\rho_1$ & Density of brine&$1.025\cdot10^{3}$& [$\mathrm{kg}\textrm{ } \mathrm{m}^{-3}$] \\
  $\frwidth $ & Fracture width &  $\EXP{\frwidth} = 5.05 \cdot 10^{-3}$ & [$\mathrm{m}$]\\ 
\hline
\end{tabular}
}
 \caption{Simulation parameters for the deterministic model problem from Section \ref{subsec_model_problem}}
  \label{table-parameter}
 \end{table}
 
%
%
\subsection{Stochastic modeling of porosity, permeability and fracture width}
\label{subsec:PorosityVar}

Stochastic modeling is used to simulate the behavior of a fractured reservoir under uncertainty. Among many possible sources of uncertainties, we consider only the following ones:  hydrogeological properties of the porous medium --- porosity ($\poro_m$) field of the surrounding porous matrix, the fracture aperture ($\frwidth$), as well as the freshwater recharge intensity $\hat{q}_x$. The permeability $K_m$ is considered to be dependent on the uncertain $\phi_m$.
The QoIs we consider are the mean and variance of the mass fraction at the six observation points.
We model the uncertain $\poro_m$ using a random field and assume $K_m$ to dependent on $\poro_m$:
\begin {eqnarray} \label {e_perm_of_poro}
 K_m = K_m (\poro_m) \in \mathbb{R}.
\end{eqnarray}
The distribution of $\poro(\bx,\xib)$, $\bx\in \D$, is determined by a set of stochastic parameters $\xib=(\xi_1,\ldots,\xi_M,...)$. Each component $\xi_i$ is a random variable.

The dependence in \eqref{e_perm_of_poro} is specific for every material.
We refer to \cite{Panda_Lake_Perm_vs_Por,Pape_Clauser_Iffland_1999,Costa_2006} for a detailed discussion. In this work, we use a Kozeny-Carman-like law
\begin{eqnarray} 
\label{e_perm_Kozeny_Carman}
 K_m (\poro_m) = \kappa_{KC} \cdot \dfrac {\poro_m^3} {1 - \poro_m^2},
\end {eqnarray}
where $\kappa_{KC}$ is a constant scaling factor. In numerical simulations we used $\kappa_{KC} = 1.5455 \cdot 10^{-8}$  [$\mathrm{m}^2$].
 
The recharge inflow flux is kept constant across the left boundary but depends on the stochastic variable $\hat{q}_{\mathrm{in}}$. We also assume that it is independent on the porosities and the permeabilities in the domain.

The uncertain width of the fracture, the recharge, and the porosity are modeled as follows, $\xi_1,\xi_2,\xi_3 \in U[-1,1]$,
\begin{equation}
\label{eq:fracture_th}
\frwidth(\xi_1) = 0.01 \cdot ((1 - 0.01) \cdot \xi_1 + (1 + 0.01)) / 2,
\end{equation}
\begin{equation}
\label{eq:recharge}
\hat{q}_{\mathrm{in}} (t, \xi_3) = 3.3 \cdot 10^{-6} \cdot (1 + 0.1 \xi_3) (1 + 0.1 \sin (\pi t / 40)),
\end{equation}
\begin{equation}
\label{eq:poro}
 \poro_m (x, y, \xi_2) = 0.35 \cdot (1 + 0.02 \cdot (\xi_2 \cos (\pi x / 2) + \xi_2 \sin (2 \pi y)).
\end{equation}
Other options for $\xib$ would be Gaussian, log-normal, beta distributions.
The corresponding permeability is defined by \eqref{e_perm_Kozeny_Carman}.

\begin{remark}
The fracture model considered would not work for open (hollow) fractures. We assume that the fracture is filled with the same porous medium regardless of its opening. This is certainly a strong assumption. But under it, the porosity and (tangential and normal) permeabilities of the fracture should be independent of its aperture. Note that the aperture is considered to be uncertain, e.g. not accurately measured but not changing as in the case of poroelasticity. However, for the model presented, the effective permeability of the entire fracture increases with aperture: Thicker fractures induce faster flow, and we observe this in our calculations.
\end{remark}

Unknown porosity can also be modeled by a random field. One of the ways to compute this random field is to make an assumption about the covariance matrix or to estimate it from the available measurement data. Then an auxiliary eigenvalue problem is solved and the truncated Karhunen-Lo\'eve expansion (KLE) is constructed \cite{Hcovariance}. Due to the high complexity of this approach, we implement a simpler one. We use $L_2$ orthogonal functions ($\sin()$ and $\cos()$), which mimic the $L_2$-orthogonal functions used in the KLE. These functions are then multiplied by the uniform random variables $\xib$ as in \eqref{eq:poro}. Of course, the actual distribution may be different. Typically, one takes the uniform or Gaussian as an initial guess and then updates it using e.g. Bayesian update techniques. But this is already the inference/inverse problem. In addition, updating the distribution function requires measurement data, which we do not have.
\section{Numerical Methods}
\label{sec:NumMethods}
\subsection{Numerical methods for the deterministic problem}
\label{sec:Num}
The system (\ref{e_cont_eq_m}--\ref{e_cont_eq_f}) is numerically solved in the domain $\D \times [0, T]$. $\D$ is covered with an unstructured grid $\D_h$ of triangles and quadrilaterals. We denote the characteristic mesh size by $h$. We apply a special technique for resolving the jump of the solution at the fracture, cf. \cite{Reiter2012} for details. The vertex-centered finite-volume scheme is used for the discretization of the system (\ref{e_cont_eq_m}--\ref{e_Darcy_vel_f}) in space (cf. \cite{Reiter2012,Reiter14}). The number of degrees of freedom associated with $\D_h$ is denoted by $n$. There are two degrees of freedom per grid vertex in $\D_h$: one for the mass fraction $\conc_{m,f}$ and another for the pressure $\pres_{m,f}$. At the fracture, several grid vertices share the same geometrical position, cf. \cite{Reiter2012,Reiter14}. We use the regular refinement to create the grid hierarchy so that $h = \mathcal{O} (n^{-1/d})$. The implicit Euler method with a time step $\tau$ is used for the time discretization of (\ref{e_cont_eq_m}--\ref{e_Darcy_vel_f}). The number of the computed time steps is $r = T / \tau$. In this paper, we focus on the mass fraction $\conc_m$ of the salt in the bulk porous medium. We denote the part of the discrete solution of the model approximating $\conc_m$ by $\sol_{m,h,\tau}$.

We use the full upwind for the convective terms (cf.\ \cite{Frolkovic-DeSchepper-ConvDom}). 
Therefore, the discretization error is of the first order w.r.t. the spatial mesh size $h$. Furthermore, the Euler method provides the first-order discretization error w.r.t. $\tau$. Thus, as $d = 2$,
\begin {equation} \label {eq:1st_ord_conv}
 \| \conc_m - \conc_{m,h,\tau} \|_2 = \mathcal{O}(h + \tau) = \mathcal{O} (n^{-1/2} + r^{-1}),
\end {equation}
which is consistent with our numerical tests.

The implicit time-stepping scheme is unconditionally stable but requires solution of the extensive nonlinear discretized algebraic system with $n$ unknowns in every time step. The Newton's method is used for this. Linear systems in the Newton iteration are solved using the BiCGStab method (cf.\ \cite {Templates}) preconditioned with the geometric multigrid method (V-cycle, cf.\ \cite{Hackbusch85}). In the multigrid cycle, the ILU${}_\beta$-smoothers \cite{Hackbusch_Iter_Sol} and Gaussian elimination are used as the coarse grid solver.
\subsection{MLMC Algorithm}
\label{sec:MLMC}

Let $\xib$ be a vector of random variables, and $g(\xib)$ the quantity of interest (QoI). In this work, we tested different QoIs: solution $\sol_m$ at a point $(t,\bx)$, integrals $I_i$ over a small subdomain $\Delta_i$ and over the whole computational domain $\D$, see Sect.~\ref{sec:numerics}. 

To estimate the propagation of uncertainties, we will calculate the mean and variance of the mass fraction. This statistical operation necessitates a large number of simulations, leading to high computational costs. The total cost can be approximated as the cost of a single simulation multiplied by the number of simulations.

To reduce the overall computational cost, we employ the well-known Multi-Level Monte Carlo (MLMC) method (see Algorithm 2). This approach is particularly suitable because the deterministic solver ug4 uses a multigrid method (see Section~\ref{sec:Num}). The MLMC method enhances efficiency by combining samples computed on different grids within the hierarchy. A more detailed discussion of these techniques can be found in~\cite{MLMC_PDE_anal11, CMLMC, giles2008, giles2015, ErikOptGeom15, teckentrup2013further, Litv_Scattered19}.

Below we provide an overview of the Multi-Level Monte Carlo (MLMC) method, drawing from various sources to ensure comprehensive coverage of implementation details. We also address aspects specific to our problem, including cost estimates at each level, convergence rates, and additional remarks. The MLMC method aims to approximate the expected value $\EXP{g}$ with an optimal computational cost. MLMC constructs a telescoping sum, defined over a sequence of spatial and temporal meshes, $\ell=0, \ldots, L$, as described next, to achieve this goal. The QoI $g$, numerically evaluated on level $\ell$, is denoted by $g_{h_{\ell},\tau_{\ell},\ell}$ or, for simplicity, by just $g_\ell$, where $h_{\ell}$ and $\tau_{\ell}$ are the discretization steps in space and time on level $\ell$.

Let $s_0$ denote the maximum computational cost of evaluating a realization of $g_0$. Similarly, let $s_\ell$ denote the computational cost of evaluating a realization of $g_\ell - g_{\ell-1}$. For simplicity, we assume that $s_\ell$ for $g_{\ell} - g_{\ell-1}$ is almost the same as $s_\ell$ for $g_{\ell}$. Since the number of iterations required is variable, the cost of computing a sample of $g_\ell - g_{\ell-1}$ may vary for different realizations.

Let $m_{\ell}$ be the number of random samples (numbers of realizations) on the $\ell$-th mesh.
The MLMC method calculates $\EXP{g_L}\approx \EXP{g}$ using the following telescopic sum: 
\begin{align}
  \EXP{g_L} &= \EXP{g_{0}} + \sum_{\ell=1}^L \EXP{g_{\ell}-g_{\ell-1}} \nonumber \\
  &\approx  
  m_0^{-1}\sum_{i=1}^{m_0} g_{0}^{(0,i)}  + \sum_{\ell=1}^L \left( m_\ell^{-1}\sum_{i=1}^{m_\ell} (g_{\ell}^{(\ell,i)} - g_{\ell-1}^{(\ell,i)} )\right). \ \label{eq:A}
\end{align}
In the equation above, the level $\ell$ in the superscript $(\ell,i)$ indicates that independent samples are used at each level. The samples used to sample $(g_{\ell} - g_{\ell-1})$ are different (independent) from the samples used to sample $(g_{\ell^{'}} -  g_{\ell^{'}-1})$, for any $\ell^{'} \neq \ell$.
The upper index $(\ell,i)$ in the notation $(g_{\ell}^{(\ell,i)} - g_{\ell-1}^{(\ell,i)})$ means that the same sample is used for $g_{\ell}$ and $g_{\ell-1}$. There are $m_0$ independent samples for $g_0$, $m_1$ independent samples for $g_{1}-g_0$, etc. Thus, in total we generated $m_0+m_1+\cdots +m_L$ independent samples.

As $\ell$ increases, the variance of $g_{\ell} - g_{\ell-1}$ decreases. Thus, the overall computational cost can be reduced by taking fewer samples on finer meshes.

\begin{remark}
To compute the second moment $\EXP{g^2}$, we proceed as follows,
$\EXP{g_L^2}\approx \EXP{g^2}$ and then to get a telescopic sum we replace $g$ by $g^2$ in \eqref{eq:A}.
\end{remark}
The information in this paragraph is problem specific. In the following numerical experiments $n_{\ell}=4n_{\ell-1}=\ldots=4^{\ell}n_{0}=2^{d\ell}n_{0}$, $d=2$, and $r_{\ell}=2r_{\ell-1}=\ldots=2^{\ell}r_{0}$. 
In the case of uniform, equidistant mesh, we could also write similar formulas for step sizes: $h_{\ell}=h_{\ell-1}\cdot 2^{-1}=h_{\ell-2}\cdot 2^{-2}=\ldots= h_0\cdot 2^{-\ell}$ and $\tau_{\ell}=\tau_0\cdot 2^{-\ell}$. The average cost $s_{\ell}$ of generating one sample of $g_{\ell}$ (the cost of one deterministic simulation for one random realization) is
\begin{equation}
s_\ell = \mathcal{O}(n_{\ell}r_{\ell})
= \mathcal{O}(4^{\ell}n_0\cdot 2^{\ell}r_0)
= \mathcal{O}(2^{2\ell}n_0\cdot 2^{\ell}r_0)
= \mathcal{O}(2^{\hat{d}\ell }n_0 r_0),
\label{eq:workpl}
\end{equation}
where $\hat{d}=d+1=3$. 

\begin{defn}
Let $Y_{\ell}:=m_{\ell}^{-1}\sum_{i=1}^{m_{\ell}} (g_{\ell}^{(\ell,i)} - g_{\ell-1}^{(\ell,i)})$, where $g_{-1}\equiv 0$, so that
\begin{align}
\label{eq:Yell} 
\EXP{Y_{\ell}}:=
      \begin{cases}
        \EXP{g_0}, & \ell=0 \\
        \EXP{g_{\ell} - g_{\ell-1}}, & \ell>0 \\
      \end{cases}.
\end{align}
Let $Y:=\sum_{\ell=0}^L Y_{\ell}$ be the multilevel estimator of $\EXP{g}$ based on $L+1$ levels and $m_{\ell}$ be the number of independent samples at level $\ell$, where $\ell=0,\ldots,L$.

Furthermore, we denote $V_0 = \Var{g_0}$ and for $\ell \ge 1$, let $V_{\ell}$ be the variance of $g_{\ell} - g_{\ell-1}$: $V_{\ell}:= \Var{g_{\ell} - g_{\ell-1}}$.
\end{defn}

The standard theory states the following facts for the mean and the variance:
\begin{equation} \label{eq:goal_variance}
 \EXP{Y}=\EXP{g_L}, \qquad \Var{Y}={\textstyle \sum_{\ell=0}^L m_{\ell}^{-1} V_{\ell}}.
\end{equation}
The cost of the multilevel estimator $Y$ is
\begin{equation}
 S := {\textstyle \sum_{\ell=0}^{L} m_{\ell}s_{\ell}}.
\end{equation}

The mean squared error (MSE) is used to measure the quality of the multilevel estimator:
\begin{equation}
\label{eq:MSE}
\MSE:=\EXP{(Y-\EXP{g})^2}=\Var{Y} + \left( \EXP{Y} - \EXP{g} \right)^2,
\end{equation}
where $Y$ is what we computed via MLMC, and $\EXP{g}$ what actually should be computed.
To achieve
$$
 \MSE \le \varepsilon^2
$$
for some prescribed tolerance $\varepsilon$, we ensure that
\begin{equation} \label{eq:bias_error}
 \left( \EXP{Y} - \EXP{g} \right)^2 = (\EXP{g_L-g})^2  \le \tfrac{1}{2} \varepsilon^2
\end{equation}
and
\begin{equation} \label{eq:var_error}
 \Var{Y} \le \tfrac{1}{2} \varepsilon^2.
\end{equation}
The bias error $\EXP{g_L-g}$ corresponds to the discretization error discussed in Sec.~ \ref{sec:Num}. Later we will see that $\EXP{Y} - \EXP{g} = O(2^{-\alpha L})$ with $\alpha \approx 1$. The bias error can be made smaller than $\varepsilon^2 / 2$ by choosing a sufficiently large $L$. Then, for this $L$, we can compute optimal $m_0, \dots, m_L$ by formula in \eqref{eq:M_ell} to provide \eqref{eq:var_error}.

In the following, we repeat the well-known \cite{giles2015} results on the computation of the sequence $m_0, \dots, m_L$. For a fixed variance $\Var{Y} =: \varepsilon^2 / 2$, the cost $S$ is minimized by choosing as $m_{\ell}$ the solution of the optimization problem

\begin{equation}
\label{eq:goal_function}
\min_{m_0,\ldots,m_{L}}F(m_0,\ldots,m_{L}),
\end{equation}
where $F(m_0,\ldots,m_{L}):=\sum_{\ell=0}^{L} \left( m_{\ell}s_{\ell}+\mu^2 \frac{V_{\ell}}{m_{\ell}}\right)$, $\mu^2$ is a Lagrange multiplier.
Thus, the derivatives of $F$ w.r.t. $m_{\ell}$ are equal to zero:
\begin{equation} \label {eq:zero_deriv}
\frac{\partial F(m_0,\ldots,m_{L})}{\partial m_{\ell}}:= s_{\ell}-\mu^2 \frac{V_{\ell}}{m_{\ell}^2}=0.
\end{equation}
Solving the system \eqref{eq:zero_deriv}, we obtain
\begin{equation} \label{eq_m_l_of_mu}
 m_{\ell}^2=\mu^2 \frac{V_{\ell}}{s_{\ell}}, \quad \text{i.e.} \quad
 m_{\ell}=\mu \sqrt{\frac{V_{\ell}}{s_{\ell}}}.
\end{equation}
Taking into account that the variation $\Var{Y}$ is fixed and substituting \eqref{eq_m_l_of_mu} into \eqref{eq:goal_variance}, i.e. $\sum_{\ell=0}^{L}V_{\ell}/m_{\ell} = \varepsilon^2 / 2$, we obtain an equation for $\mu$:
$$
\sum_{\ell=0}^{L}\frac{V_{\ell}}{\mu \sqrt{\frac{V_{\ell}}{s_{\ell}}}} = \tfrac{1}{2} \varepsilon^2.
$$
From this equation, we get $\mu = 2 \varepsilon^{-2} \sum_{\ell=0}^{L} \sqrt{V_{\ell}s_{\ell}}$, and therefore
\begin{equation}
\label{eq:M_ell}
 m_{\ell} = 2 \varepsilon^{-2} \cdot \sqrt{\frac{V_{\ell}}{s_{\ell}}} \cdot \sum_{i=0}^{L} \sqrt{V_{i}s_{i}}.
\end{equation}
For this set of $m_{\ell}$, the total computational cost of $Y$ is
\begin{equation}
\label{eq:total_cost_MLMC}
S = 2 \varepsilon^{-2} \left( \sum_{\ell=0}^L \sqrt{V_{\ell} s_{\ell}}\right)^2.
\end{equation}
For further analysis of this sum, see \cite{giles2015}, p.4.


The cost increases exponentially with $\ell$ while the weak error $\EXP{g_L-g}$ and multilevel correction variance $V_{\ell}$ decrease exponentially leads to the following theorem (cf. Theorem 1, p.~6 in \cite{giles2015}):
\begin{theorem}
\label{thm:costMLMC}
Consider a fixed $t=t^*$. Suppose positive constants $\alpha,\beta,\gamma > 0$ exist such that $\alpha \geq \frac{1}{2} \min(\beta, \gamma \hat{d})$, and
\begin{subequations}
\label{eq:q1q2}
\begin{align}
    \vert \EXP{g_\ell-g} \vert & \leq c_1 2^{-\alpha\ell} \label{eq:weak_error_model} \\
    V_{\ell} & \leq c_2 2^{-\beta\ell} \label{eq:strong_error_model} \\
     s_{\ell} &\leq c_3 2^{\hat{d}\gamma \ell}.
\end{align}
\end{subequations}
Then, for any accuracy $\varepsilon < e^{-1}$, a constant $c_4>0$ and a sequence of realizations $\{m_{\ell}\}_{\ell=0}^L$ exist, such that
$\MSE < \varepsilon^2$, where $\MSE$ is defined in \eqref{eq:MSE},
%
and the computational cost is
\begin{align}
\label{eq:mlmc_iso_work} 
S=
      \begin{cases}
        c_4 {\varepsilon^{-2}}, & \beta > \hat{d}\gamma \\
        c_4 {\varepsilon^{-2} \left(\log(\varepsilon)\right)^2}, & \beta= \hat{d}\gamma \\
        c_4 {\varepsilon^{-\left(2 +\frac{\hat{d}\gamma-\beta}{\alpha}\right)}},  & \beta < \hat{d}\gamma. \\
      \end{cases}
\end{align}
\end{theorem}
The parameter $\gamma$ is the complexity rate of the numerical solver. Our geometric multigrid solver has linear complexity, i.e. $\gamma=1$.

This theorem (see also \cite{hoel2014implementation,hoel2012adaptive,charrier2013,MLMC_PDE_anal11,giles2008}) indicates that, even in the worst-case scenario, the MLMC algorithm has a lower computational cost than that of the traditional (single-level) MC method, which scales as $\mathcal{O}(\varepsilon ^{-2-\hat{d}\gamma/\alpha})$. 
\begin{remark}
In Theorem~\ref{thm:costMLMC}, the factors $c_1$, $c_2$, $c_3$, $c_4$ as well as the exponents $\alpha$, $\beta$ and $\gamma$  depend on the time point $t$. This makes $L$ and $m_\ell$ time-dependent, too.
\end{remark}
\begin{remark}
One possible choice for a scaling factor in \eqref{eq:bias_error} and \eqref{eq:var_error} is $E_0:=|\EXP{g_0(t^*, \bx^*)}|$, where $g_0(t^*, \bx^*)$ is the solution computed on level $\ell=0$ at point $(t^*, \bx^*)$.
\end{remark}
\begin{remark}
We consider the error relatively to $g$. For this, in \eqref{eq:bias_error} and \eqref{eq:var_error}, we replace $\varepsilon$ by $\varepsilon \cdot E_0$. Equivalently, we can divide (rescale) $\EXP{g_L-g}$ by $E_0$ and $V_i$ by $E_0^2$. Therefore, we get:
\begin{equation}
\label{eq:Egl}
{ | \EXP{g_L-g} |} \leq c_1 2^{-\alpha L}.
\end{equation}
Now, to satisfy \eqref{eq:bias_error}, we want 
\begin{equation}
\label{eq:rel_mean_eps}
{| \EXP{g_L-g} |}\leq  \frac{1}{\sqrt{2}} \varepsilon E_0.
\end{equation}
From this inequality, we can estimate $L$:
\begin{equation}
\label{eq:est_L}
c_1 2^{-\alpha L} = \frac{\varepsilon E_0  }{\sqrt{2}}
\end{equation}
\begin{equation} \label{eq:num_levels}
 L = - \tfrac{1}{\alpha} \log_2 \frac {\varepsilon E_0 } {\sqrt{2} C_1}.
\end{equation}
Equations \eqref{eq:M_ell} and \eqref{eq:total_cost_MLMC} attain the form
\begin{equation} \label{eq:scaled_m_l_S}
 m_{\ell} = \frac {2 \varepsilon^{-2}} {{E_0} ^2} \cdot  \sqrt{\frac{V_{\ell}}{s_{\ell}}} \cdot \sum_{i=0}^{L} \sqrt{V_{i}s_{i}},
 \qquad
 S = \frac {2 \varepsilon^{-2}} {E_0^2} \left( \sum_{\ell=0}^L \sqrt{V_{\ell} s_{\ell}}\right)^2.
\end{equation}
\end{remark}

Using preliminary numerical tests (see Fig.~\ref{fig:weak_strong_p1_t15}) and following preprocessing steps as in Algorithm~\ref{alg:MLMC_prep}, we can estimate the convergence rates $\alpha$ for the mean (the so-called weak convergence) and $\beta$ for the variance (the so-called strong convergence), as well as the constants $c_1$ and $c_2$ from \eqref{eq:weak_error_model}-\eqref{eq:strong_error_model}.
In addition, the rate $\alpha$ is strongly related to the order of the discretisation error (see Section~\ref{sec:Num}), which is equal to $1$. Note that precise estimates of the parameters $\alpha$ and $\beta$ are crucial for the optimal complexity of the MLMC method and for the optimal distribution of the computational effort over $L+1$ layers. The number of samples $m_{\ell}$ and the level $L$ will be suboptimal if the rates $\alpha$ and $\beta$ are not accurately estimated. To obtain these estimates, we run the Algorithm~\ref{alg:MLMC_prep}. Some advanced implementations of MLMC, e.g. CMLMC \cite{CMLMC}, can update the convergence rate on the fly. The calculated rates are used in the Theorem~\ref{thm:costMLMC} and the Algorithm~\ref{alg:MLMC}. This Algorithm~\ref{alg:MLMC} estimates the mean value of the QoI. Examples of calculated and used parameters $m_{\ell}$ and $L$ are listed in the rows of Table~\ref{tab:ml}. The values in Table~\ref{tab:ml} are only valid for a fixed QoI. Both algorithms should be rerun, and all estimates made again for a new QoI.

\begin{algorithm}
\caption{Preprocessing for MLMC algorithm}\label{alg:MLMC_prep}
\begin{algorithmic}
\State \textbf{Input:} set, for example, $L=3$, $m_{\ell}=10$, $\ell=0,\ldots,L$ 
\State Drawn $g_{\ell}(\omega_i) - g_{\ell - 1}(\omega_i)$ for $i=1\ldots m_{\ell}$ and $\ell=0\ldots L$
\State Estimate convergence rate of $\EXP{g_{\ell}(\omega_i) - g_{\ell - 1}(\omega_i)}$, $\Var{g_{\ell}(\omega_i) - g_{\ell - 1}(\omega_i)}$ for each $\ell$
\State \textbf{Output:} Averaged estimated convergence rates $\alpha$ and $\beta$, complexity rate $\gamma$, constants $c_1,c_2,c_3$
\end{algorithmic}
\end{algorithm}
\begin{algorithm}
\caption{MLMC algorithm}\label{alg:MLMC}
\begin{algorithmic}
\State \textbf{Input:} MSE Error $\varepsilon^2$
\State Estimate $L$ as in \eqref{eq:num_levels}
\State Compute $m_{\ell}$, $\ell=0,\ldots,L$ as in \eqref{eq:scaled_m_l_S}
\State Compute $\EXP{g_{\ell}(\omega_i) - g_{\ell - 1}(\omega_i)}$, for every $\ell$ and $i=1\ldots m_{\ell}$ 
\State Build the telescopic sum as in \eqref{eq:A}.
\State \textbf{Output:} the estimated mean value.
\end{algorithmic}
\end{algorithm}
\section{Numerical Experiments}
\label{sec:numerics}

The simulations below depend on three stochastic variables $\{\xi_1,\xi_2,\xi_3\}$, 
$\xi_1$ is controlling thickness of the fracture \eqref{eq:fracture_th},
$\xi_2$ the values of permeability and porosity \eqref{eq:poro}, and
$\xi_3$ the intensity of the recharge, as described in \eqref{eq:recharge}.

Thus, in total there are four uncertain parameters: thickness, recharge, permeability and porosity. Permeability and porosity are mutually dependent. To verify that the uncertainties in all four parameters contribute significantly to the solution, we performed auxiliary experiments. These experiments showed that all parameters contribute significantly and that there are no parameters whose contribution could be neglected.


We consider the following QoIs: 1) the solution $\conc_m$ at a point $(t,\bx)$, i.e. $g=\conc_m(t,\bx)$, and 2) an integral over a small sub-domain 
\begin{equation}
\label{eq:integral_box}
I_i(t,\omega):=\int_{\bx\in \Delta_i} \conc_m(t,\bx,\omega) \rho (\conc_m(t,\bx,\omega))  d\bx,\quad \text{where}
\end{equation}
\begin{equation}
\Delta_i := [x_i-0.1,x_i+0.1]\times [y_i-0.1,y_i+0.1],\quad i=1,\ldots,6.
\end{equation}
The list of all points $\bx_i=(x_i,y_i)$ is defined in \eqref{eq:6points}.
The value of $I_i$ is the mass of the salt in a subdomain $\Delta_i$. The size of each $\Delta_i$ is small ($0.2^2=0.04$), compared to $\D$.

The numerical scheme provides only the first order of the accuracy to compute $\conc_m(t,\bx_i,\omega)$ and $I_i(t,\omega)$, i.e., the convergence rate $\alpha$ (weak convergence in \eqref{eq:weak_error_model}) should be $\approx 1$. The numerical results below show that, indeed, $\alpha \approx 1$.

Table~\ref{tab:comp_times_2params_fracture} contains computing times needed to compute the solution $\sol_m$ at each level $\ell$. The fifth column contains the average computing time, and the sixth and seventh columns contain the shortest and longest computing times. The computing time for each simulation varies depending on the number of iterations, which depends on the porosity and permeability. We observed that, after $\approx 6016$~sec., the solution is almost unchanging; thus, we perform the experiment only for $t\in [0, T]$, where $T=6016$ sec. For example, if the number of time steps is $r_{\ell}=188$ (Level 0 in Table~\ref{tab:comp_times_2params_fracture}), then the time step $\tau = \frac{T}{r_{\ell}}=\frac{6016}{188}=32$~sec. 

\begin{table}[htbp!]
\centering
\begin{tabular}{|l|l|l|l|l|l|l|}
\hline 
\multirow{2}{*}{Level $\ell$}& 
\multirow{2}{*}{$n_{\ell}$, ($\frac{n_{\ell}}{n_{\ell-1}}$)} &
\multirow{2}{*}{$r_{\ell}$, ($\frac{r_{\ell}}{r_{\ell-1}}$)} &
\multirow{2}{*}{$\tau_{\ell}=6016/r_{\ell}$} &
\multicolumn{3}{|c|}{Computing times ($s_{\ell}$), ($\frac{s_{\ell}}{s_{\ell-1}}$)}\\
    &&   &                               & average & min. & max.\\
\hline
0 &    608      &   188     & 32 &     3      &    2.4 &     3.4 \\ \hline 
1 &   2368 (3.9)  &  376 (2) & 16 &     22 (7.3) &    15.5 &     27.8 \\ \hline 
2 &  9344 (3.9)& 752 (2) &  8 &   189 (8.6) &  115 &   237 \\ \hline 
3 &  37120 (4)& 1504 (2) &  4 & 1831 (10) & 882 & 2363 \\ \hline 
4 & 147968 (4) & 3008 (2) &  2 & 18580 (10) & 7865 & 25418 \\ \hline 
\end{tabular}
\captionsetup{width=.99\textwidth}
\caption{Number of degrees of freedom $n_{\ell}$, number of time steps $r_{\ell}$, step size in time $\tau_{\ell}$, average, minimal, and maximal computing times on each level $\ell$.} 
\label{tab:comp_times_2params_fracture} 
\end{table}

\begin{figure}[htbp!]
\begin{center}
\includegraphics[width=0.5\textwidth]{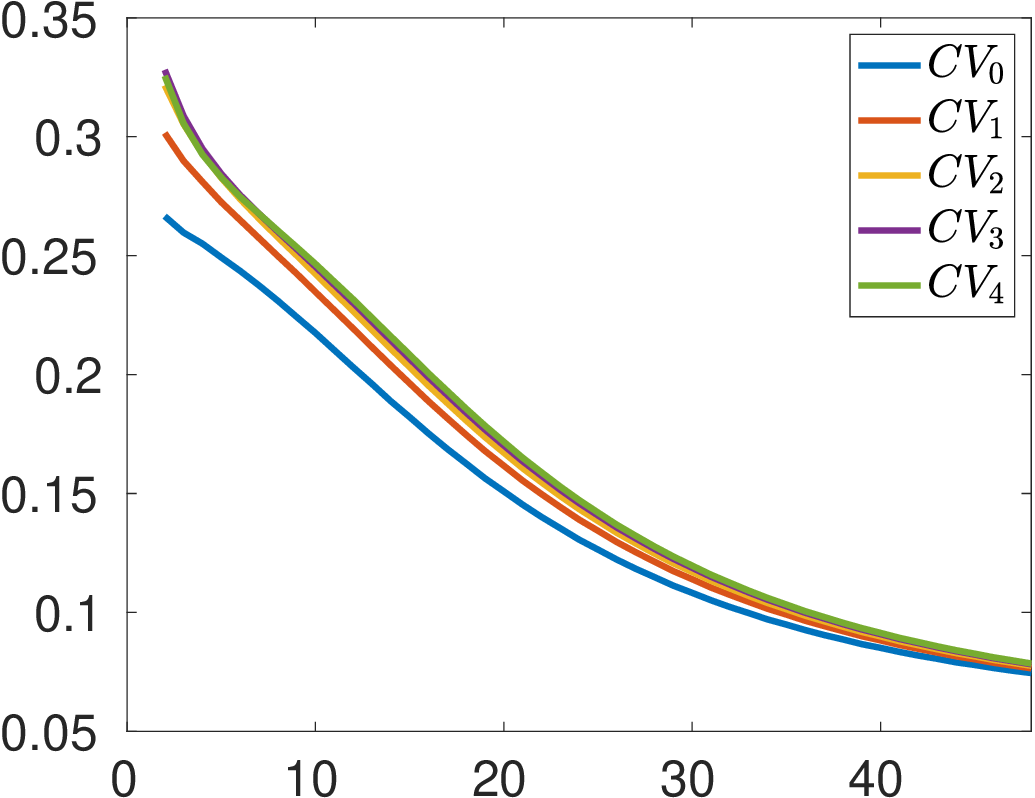}
\includegraphics[width=0.47\textwidth]{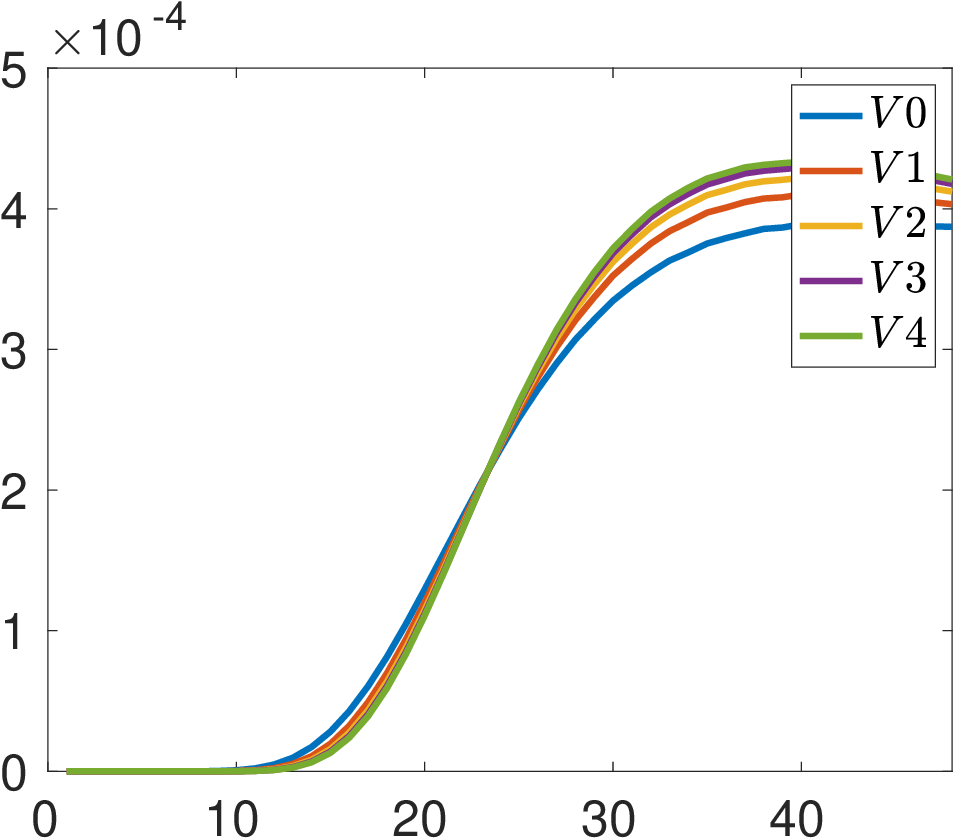}\\
\includegraphics[width=0.49\textwidth]{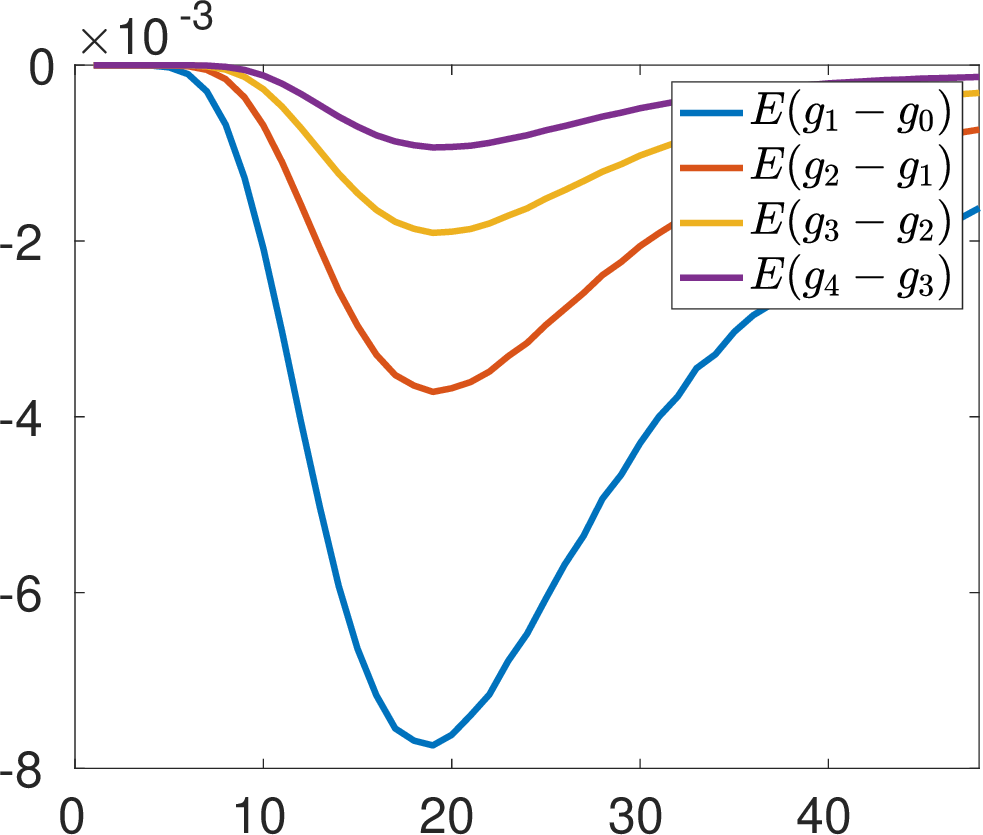}
\includegraphics[width=0.49\textwidth]{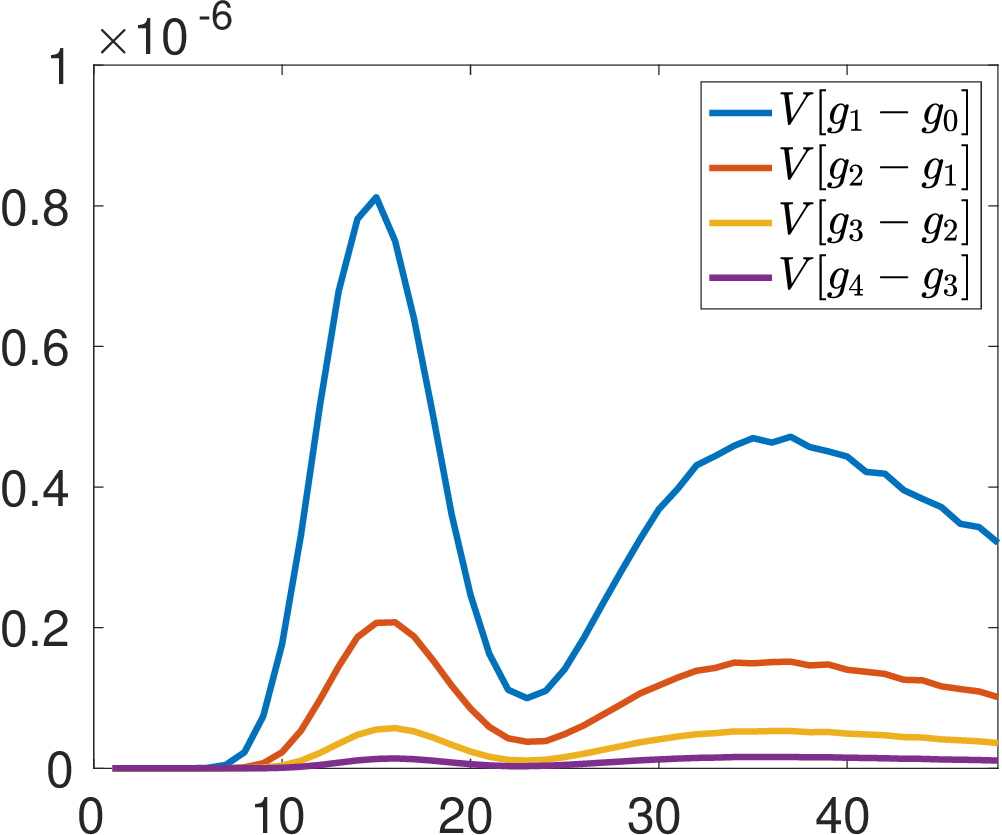}
   \caption{(top left) The coefficient of variance $CV_{\ell}:=CV(g_{\ell})$, (top right) the variance $\Var{g_{\ell}}$, (bottom left) the mean $\EXP{g_{\ell} - g_{\ell-1}}$, (bottom right) the variance $\Var{g_{\ell} - g_{\ell-1}}$. The QoI is $g=\sol_m(t,\bx_1)$. The small oscillations in the two lower pictures are due to the dependence of the recharge $\hat{q}_{\mathrm{in}}$ on the time, cf. \eqref{eq:recharge}. The time is along the $x$ axis, $t\in[0,47\tau]$}.
    \label{fig:QOIs_for_point1}
\end{center}
\end{figure}

In Figure~\ref{fig:QOIs_for_point1} (top left) we visualise the coefficient of variation $CV_{\ell}(g):=\frac{\sigma(g_{\ell})}{\EXP{g_{\ell}}}$, which is the ratio of the standard deviation $\sigma(g_{\ell})$ to the mean $\EXP{g_{\ell}}$. The variance $V_{\ell}:=\Var{g_{\ell}}$ is shown on the top right. The mean value of the difference $\EXP{g_{\ell} - g_{\ell-1}}$ is shown on the bottom left and the variance $\Var{g_{\ell} - g_{\ell-1}}$ on the bottom right. The QoI is $g:=\sol_m(t, \bx_1)$ and the time $t=0,\ldots,47\tau$ is along the $x$ axis.

\subsection{The mean value and variance}

\begin{figure}[t!]
\begin{center}
\includegraphics[width=0.24\textwidth]{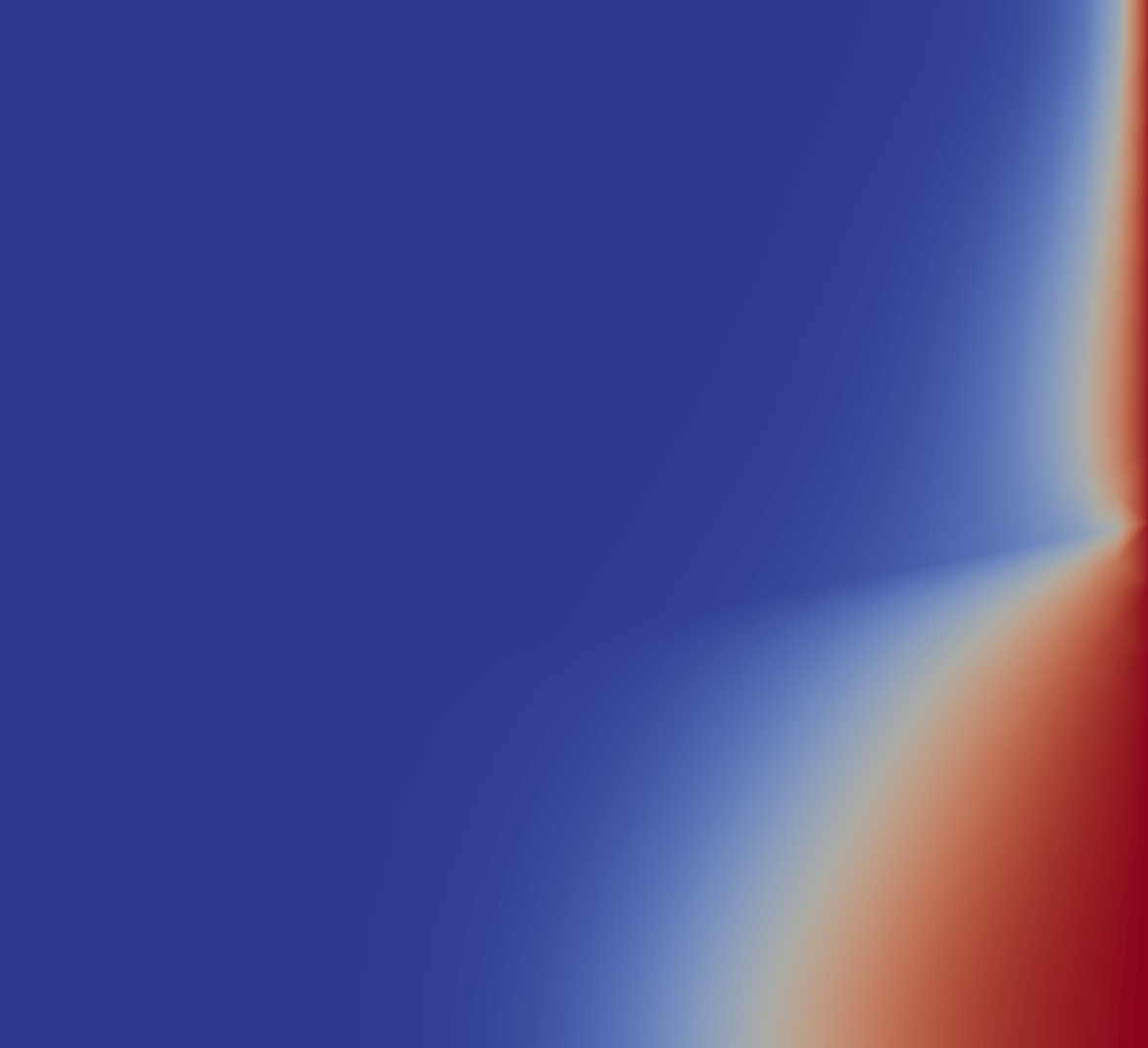}
\includegraphics[width=0.24\textwidth]{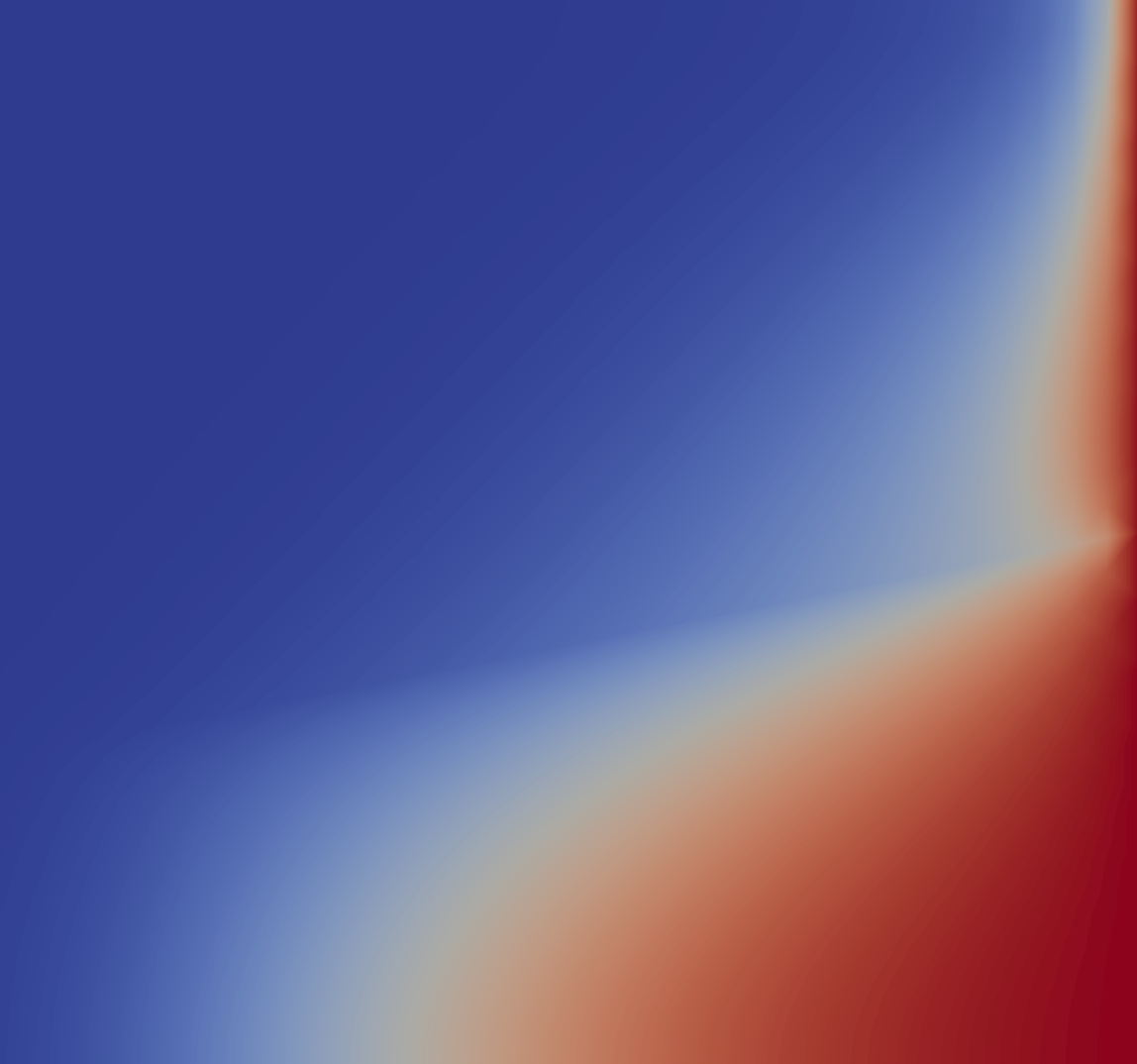}
\includegraphics[width=0.24\textwidth]{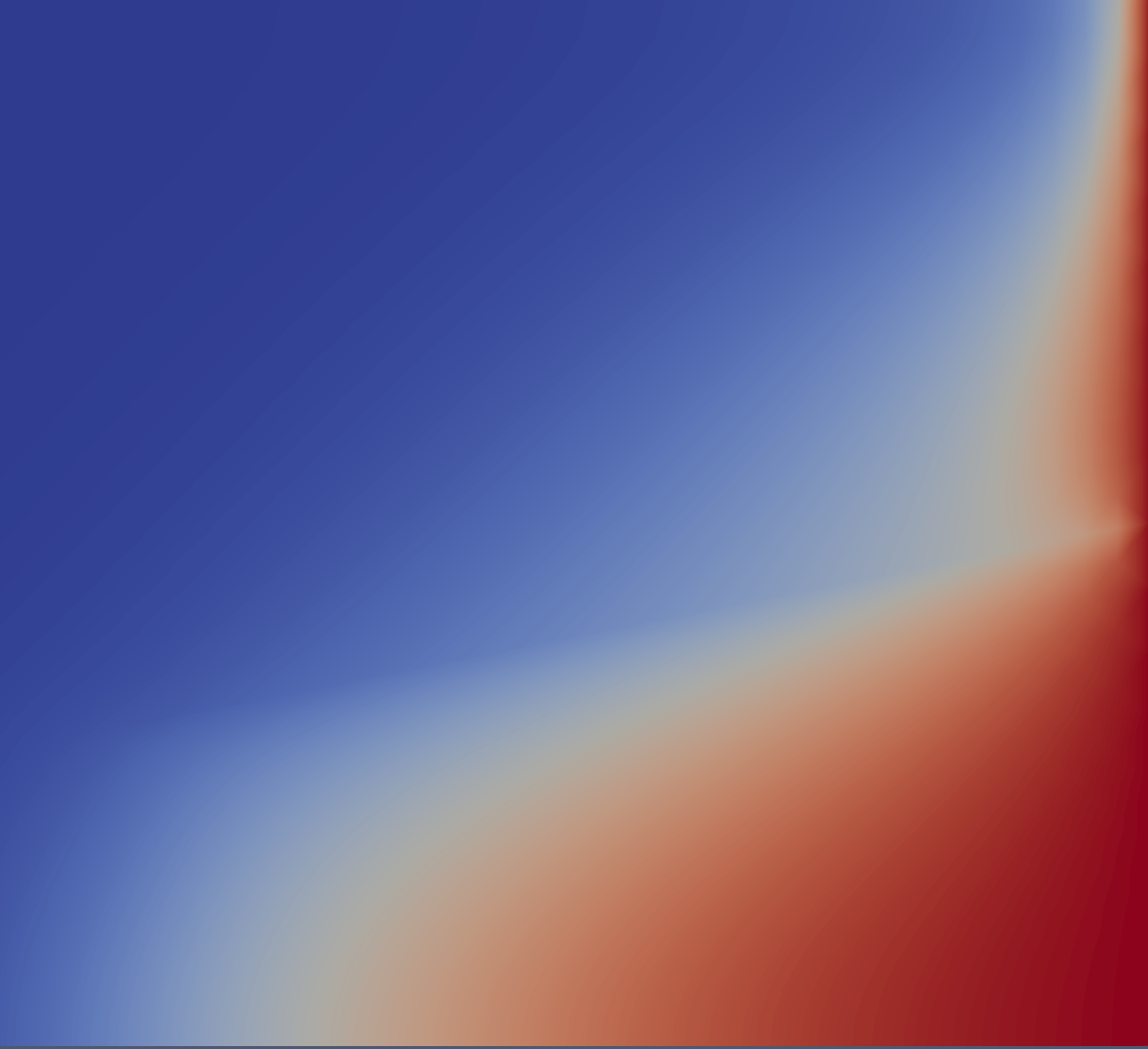}
\includegraphics[width=0.24\textwidth]{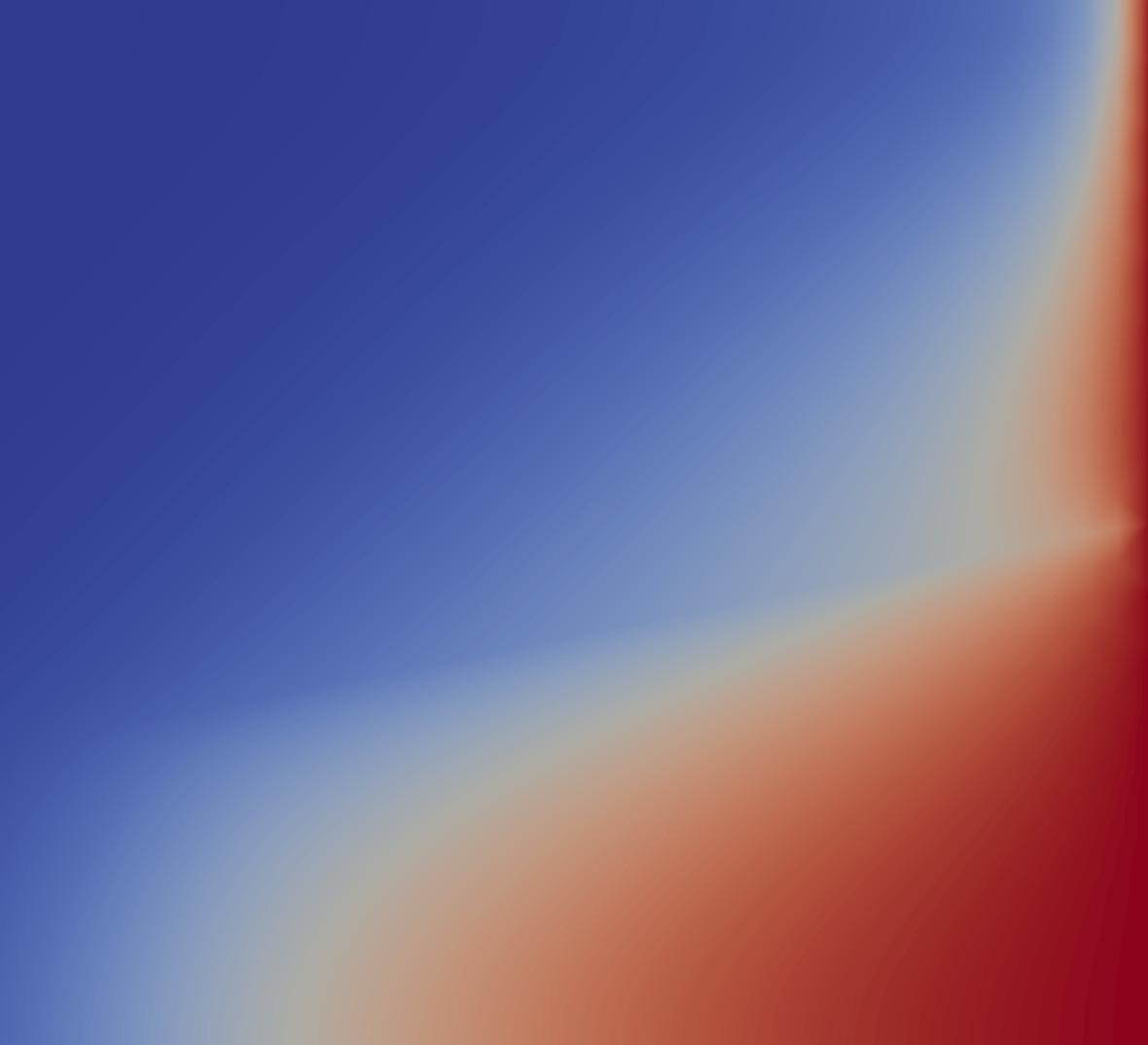}
\caption {The mean value $\EXP{\sol_m(t,\bx)}$ at $t=\{7\tau, 19\tau, 40\tau, 94\tau\}$. In all cases, $\EXP{\sol_m}(t,\bx) \in [0,1]$. The blue colour corresponds to $\EXP{\sol_m}(t,\bx)=0$ and
the dark red colour to $\EXP{\sol_m}(t,\bx)=1$.}
\label{fig:mean_T}
\end{center}    
\end{figure}

Figure~\ref{fig:mean_T} shows the mean value $\EXP{\sol_m(t,\bx)}$ at time points $\{7\tau, 19\tau, 40\tau, 94\tau\} $. In all cases $\EXP{\sol_m(t,\bx)} \in [0,1]$, the dark red color corresponds to the value 1 and dark blue to the value 0. The fracture is visible in the center, too. Due to the high permeability inside the fracture for the most realizations, the flow is divided into two subdomains (compare the red areas above and below the fracture). The high-permeable fracture isolates the upper flow from the lower flow. If we further increase the permeability inside the fracture, we get two almost independent problems: one above the fracture and one below.

\begin{figure}[t!]
\begin{center}
\includegraphics[width=0.24\textwidth]{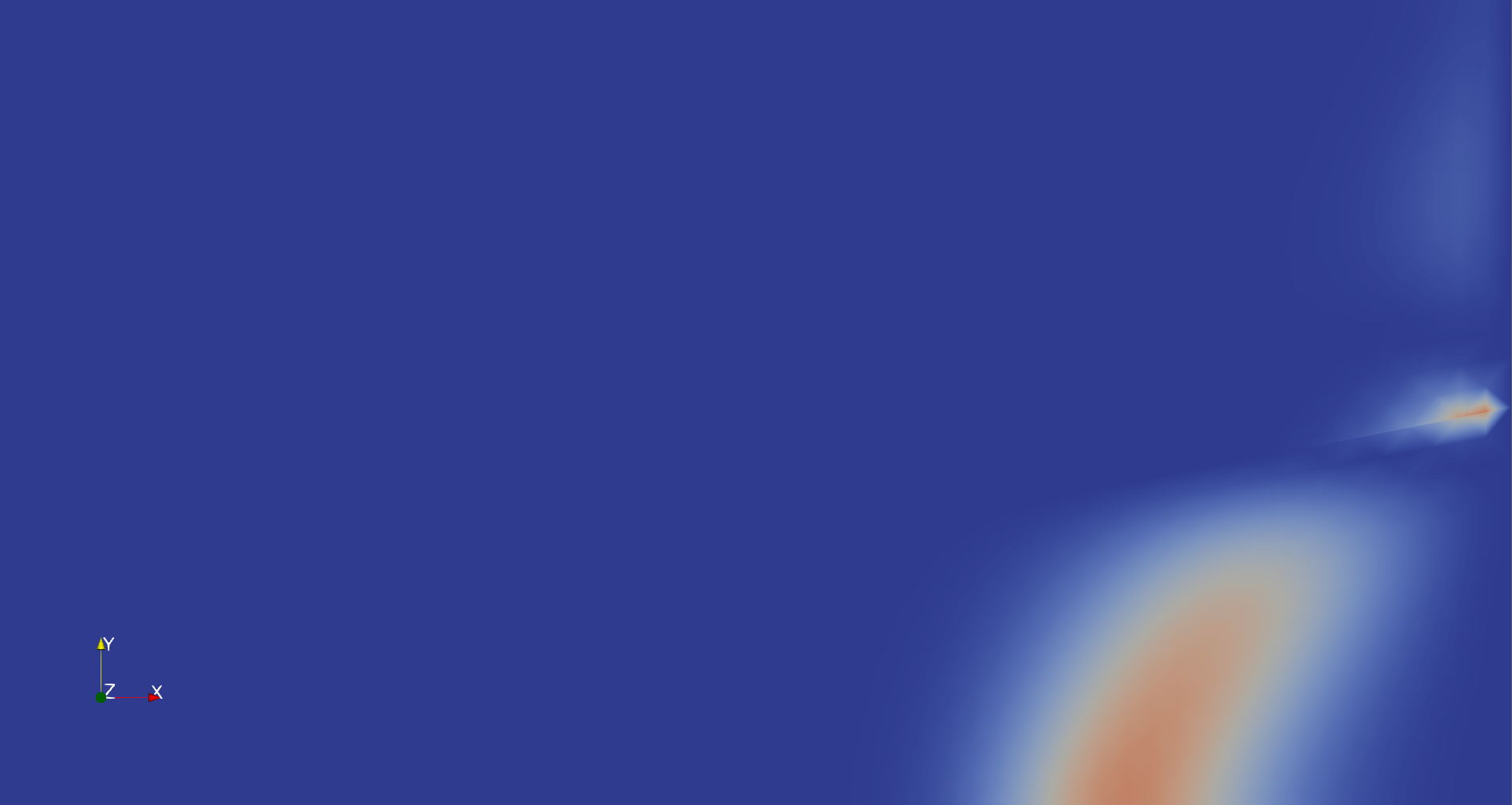}
\includegraphics[width=0.24\textwidth]{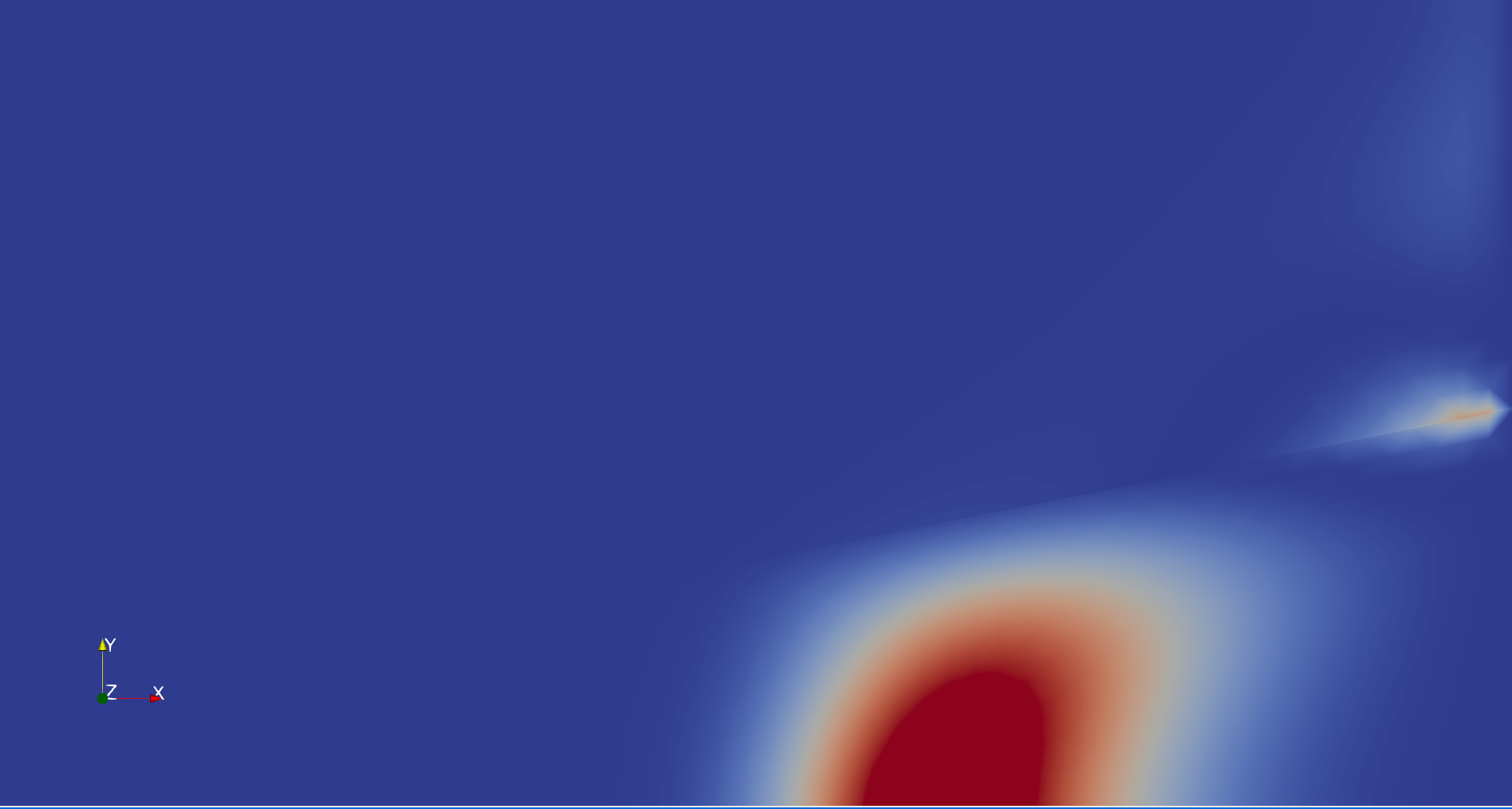}
\includegraphics[width=0.24\textwidth]{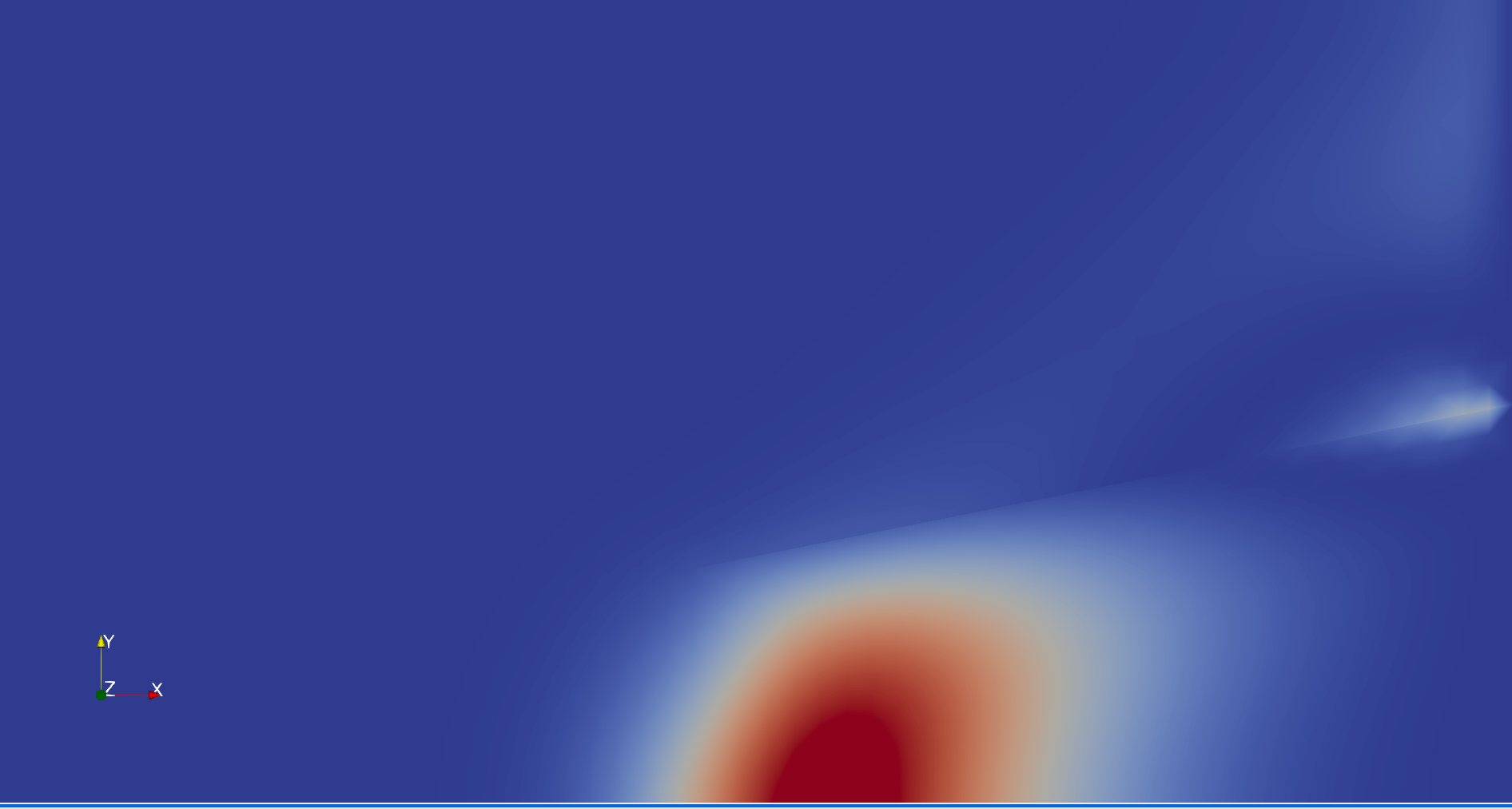}
\includegraphics[width=0.24\textwidth]{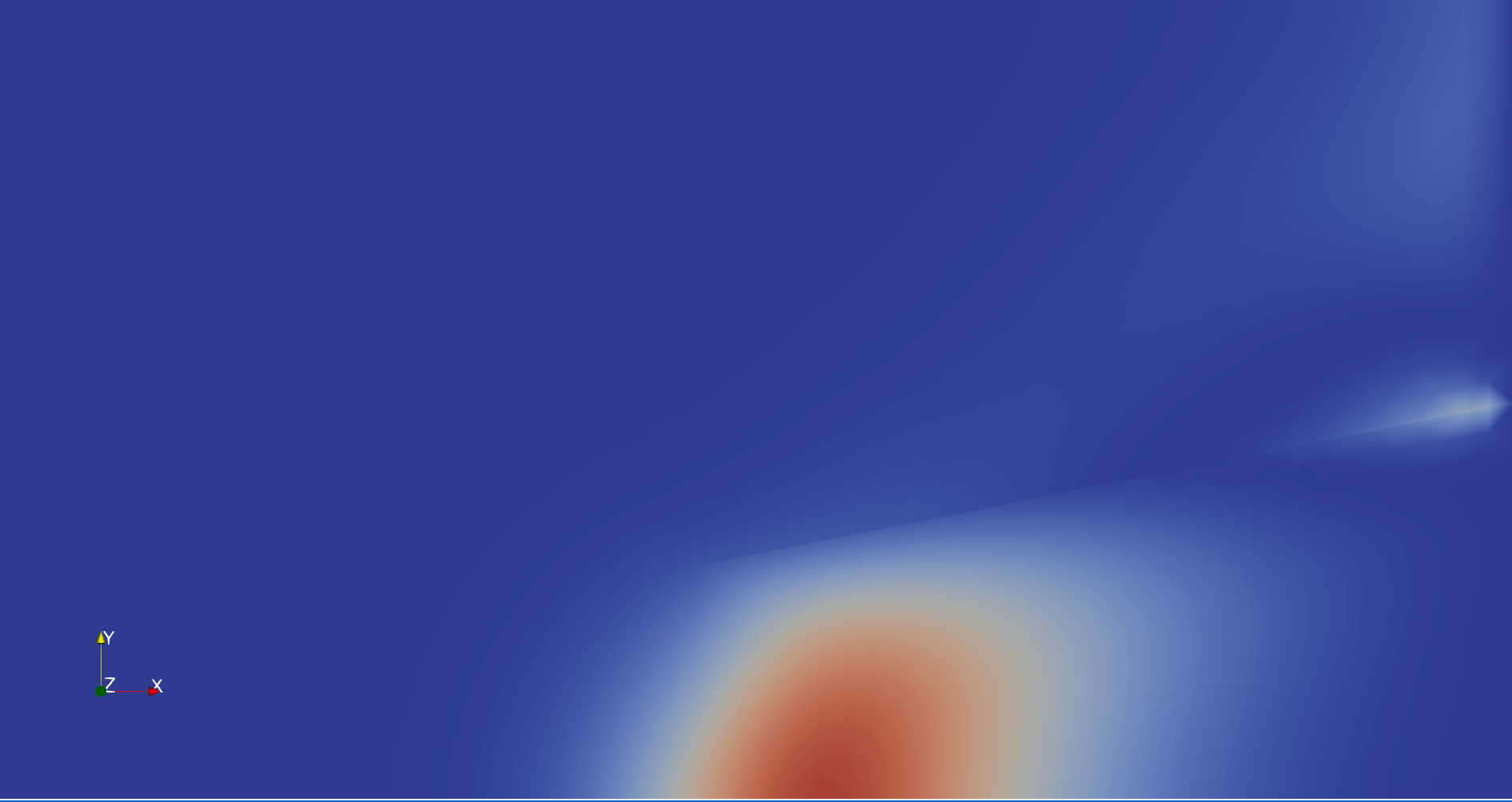}
\caption {The variance $\Var{\sol_m(t,\bx)}$ at $t=\{7\tau, 19\tau, 40\tau, 94\tau\}$. Maximal values (dark red colour) of $\Var{\sol_m}$ are $1.9\cdot 10^{-3}$, $3.4\cdot 10^{-3}$, $2.9\cdot 10^{-3}$, $2.4\cdot 10^{-3}$ respectively. The blue colour corresponds to a zero value.}
\label{fig:var_T}
\end{center}    
\end{figure}

Similarly, Fig.~\ref{fig:var_T} presents the variance of the salt mass fraction in the computational domain for $t=\{7\tau, 19\tau, 40\tau, 94\tau\}$. The maximal values of $\Var{\sol_m}$ (denoted by the dark red colour) are $1.9\cdot 10^{-3}$, $3.4\cdot 10^{-3}$, $2.9\cdot 10^{-3}$, $2.4\cdot 10^{-3}$ respectively.

The general flow and transport patters in this setting were discussed in Sect.~\ref {subsec_model_problem}. The aperture $\epsilon$ of the fracture has the most essential influence on the $\sol_m$. If $\epsilon$ is small (i.e.\ for $\xi_1 = -1$, cf. \eqref{eq:fracture_th}), the fracture does not increase the permeability of the medium significantly so the flow and transport are very similar to the Henry problem. For thick fractures $(\xi_1 = 1)$, the separation of the upper and the lower part of the domain (see Sect.~\ref {subsec_model_problem}) is very strong. This is clearly visible in the mean value of $\sol_m(t,\bx)$ in Fig.~\ref{fig:mean_T}. Note that $\epsilon$ changes the flow under the left end of the fracture where fluid enters it. This area is characterized by the high value of the variance, see Fig.~\ref{fig:var_T}. The influence of the other uncertainties is weaker but nevertheless important, too.

\newpage
\subsection{Estimation of convergence rates $\alpha$ and $\beta$ for various QoIs}
\begin{figure}[t!]
\begin{center}
\includegraphics[width=0.49\textwidth]{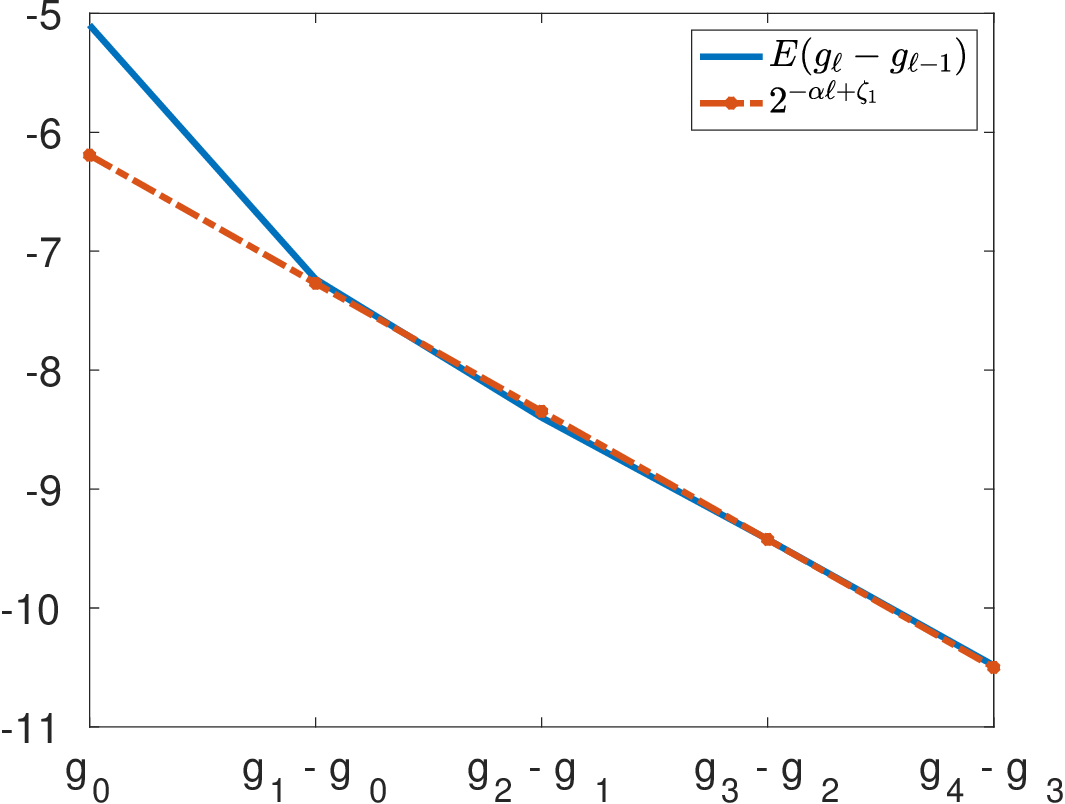}
\includegraphics[width=0.49\textwidth]{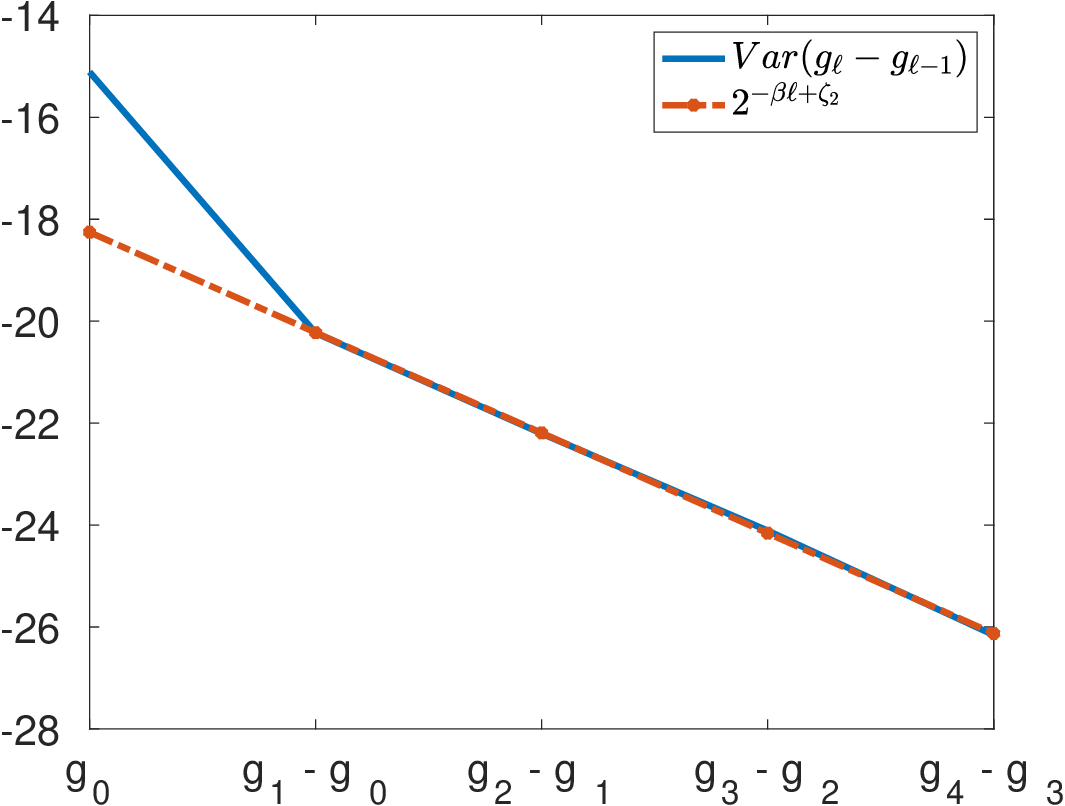}
\caption {The weak and the strong convergences, QoI is $g:=\sol_m(t_{15},\bx_1)$, $\alpha=1.07$, $\zeta_1=-1.1$, $\beta=1.97$, $\zeta_2 = -8$. The $0y$ axis denotes $\log_2(\EXP{g_{\ell}- g_{\ell-1}})$ on the left and $\log_2(\Var{g_{\ell}- g_{\ell-1}})$ on the right.}
\label{fig:weak_strong_p1_t15}
\end{center}    
\end{figure}

Figure~\ref{fig:weak_strong_p1_t15} (left) shows the decay of $\EXP{g_{\ell} - g_{\ell-1}}$ as a function of $g_{\ell} - g_{\ell-1}$. The $0y$ axis denotes $\log_2(\EXP{g_{\ell}- g_{\ell-1}})$. This dependence is fitted by a $2^{-\alpha \ell + \zeta_1}$ function. The parameter $\alpha$ indicates the weak convergence rate. Figure~\ref{fig:weak_strong_p1_t15} (right) shows the decay of $\Var{g_{\ell} - g_{\ell-1}}$ with respect to $g_{\ell} - g_{\ell-1}$. The $0y$ axis denotes $\log_2(\Var{g_{\ell}- g_{\ell-1}})$.
This dependence is fitted by a $2^{-\beta \ell + \zeta_2}$ curve, where the parameter $\beta$ indicates the strong convergence rate. The constants are $C_1=0.47$ and $C_2=3.7\cdot 10^{-3}$. The QoI is $g:=\sol_m(t_{15},\bx_1)$. The calculated rates $\alpha=1.07$ and $\beta=1.97$ are very close to the theoretical values, which are $1$ and $2$ respectively.

\begin{figure}[t!]
\begin{center}
\includegraphics[width=0.49\textwidth]
{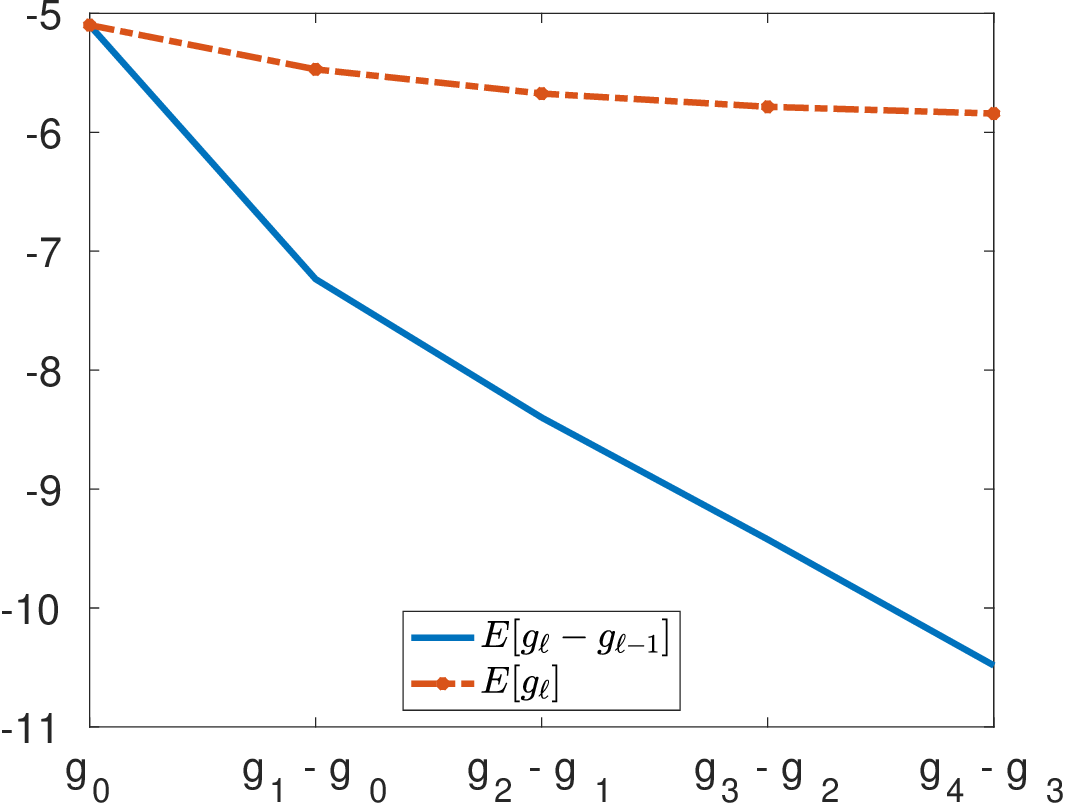}
\includegraphics[width=0.49\textwidth]
{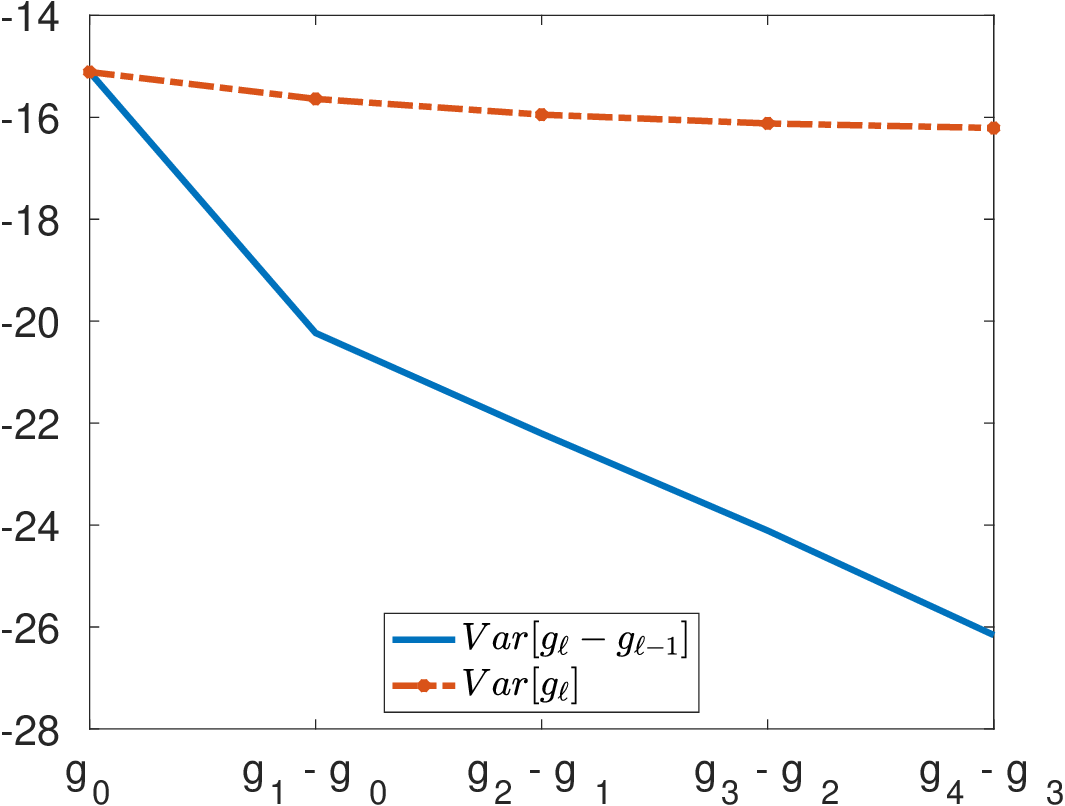}
\caption{Decay comparison of (left) $\EXP{g_{\ell} - g_{\ell-1}}$ and $\EXP{g_{\ell}}$ vs. $\ell$; (right) $\Var{g_{\ell} - g_{\ell-1}}$ and $\Var{g_{\ell}}$. QoI $g= \sol_m (t_{15},\bx_1)$, $\bx_1=(1.1, -0.8)$. The $0y$ axis shows $\log_2(\cdot)$ of $\EXP{g_{\ell}- g_{\ell-1}}$, $\EXP{g_{\ell}}$, $\Var{g_{\ell}- g_{\ell-1}}$, and $\Var{g_{\ell}}$ values.}
\label{fig:MLMC_MC_mean}
\end{center}    
\end{figure}

In Fig.~\ref{fig:MLMC_MC_mean} (left) the decays of $\EXP{g_{\ell} - g_{\ell-1}}$ and $\EXP{g_{\ell}}$ are plotted against the level $\ell$. The QoI is $\sol_m(t_{15},\bx_1)$, $\bx_1=(1.1, -0.8)$. The values on the $0y$ axis are $\log_2(\EXP{g_{\ell}- g_{\ell-1}})$ and  $\log_2(\EXP{g_{\ell}})$.
Similarly, Figure~\ref{fig:MLMC_MC_mean} (right) shows that $\Var{g_{\ell} - g_{\ell-1}}$ decays much faster than $\Var{g_{\ell}}$. The values on the $0y$ axis are $\log_2(\Var{g_{\ell}- g_{\ell-1}})$ and  $\log_2(\Var{g_{\ell}})$.

In Table~\ref{tab:estimated_rates} we provide estimations for the weak and strong convergence rates $\alpha$ and $\beta$, as well as constants $C_1$ and $C_2$ for different QoIs. Each QoI is the solution at a point $\bx_i$ or an integral $I_i$ as in \eqref{eq:integral_box}. We see that the value $\alpha$ for $\sol(\bx_i, 15\tau)$ is very close to $-1.0$, and $-1.5$ for $I_2$. The value $\alpha=1.0$ is fully corresponding to the theory. The estimated value of $\beta$ is also very close to the theoretical value, which is $2\alpha\approx 2$. The values of $\alpha$ and $C_1$ are used for estimation of the number of needed levels $L$, see also \cite{LOGASHENKO2024_JCP}. The value $\alpha=-1.5$ for $I_2$ is also what we expected, since the integral of $\sol$ is usually a smoother function than just $\sol$. The value $\beta=-2.6$ (in the last row) is also close to the theoretical value, which is $2\alpha = -3$.

\begin{table}[htbp!]
\centering
\begin{tabular}{|c|c|c|c|c|} \hline
  QoI & $\alpha$ & $c_1$ & $\beta$ & $C_2$ \\ \hline
  $\sol_m(\bx_1, 15\tau)$ & -1.07 & 0.47& -1.97& $4\cdot 10^{-3}$ \\ \hline
  
  $\sol_m(\bx_2, 15\tau)$ &-0.96& 0.25& -1.6& $3.1\cdot 10^{-4} $ \\ \hline
  
  $\sol_m(\bx_3,15\tau)$ & -0.9& 0.12& -1.9&  $4.5\cdot 10^{-5} $ \\ \hline
  
  $\sol_m(\bx_4,15\tau)$ & -1.07& 0.47& -1.94& $3.6\cdot 10^{-3} $  \\ \hline
  
  $\sol_m(\bx_5,15\tau)$ & -0.95& 0.25& -1.5& $1.9\cdot 10^{-4} $  \\ \hline
  
  $\sol_m(\bx_6,15\tau)$ & -0.89& 0.1& -1.82& $9.2\cdot 10^{-5} $  \\ \hline
  
$I_2(\sol_m)$ & -1.5 & 8.3& -2.6& 0.4 \\ \hline
    \end{tabular}
\captionsetup{width=.99\textwidth}
\caption{Estimated the weak and strong convergence rates $\alpha$ and $\beta$, and constants $C_1$ and $C_2$ for different QoIs. QoI are solution in a point $\bx_i$ and an integral $I_2$ as in \eqref{eq:integral_box}.}
 \label{tab:estimated_rates}
\end{table}

%
In Table~\ref{tab:ml} we show the numbers $m_{\ell}$ needed to get the MSE error below $\varepsilon^2$, and $\varepsilon$ should be understood as a relative error, as it is defined in \eqref{eq:rel_mean_eps}. All numbers are obtained for the QoI $g=\sol_m(t_{15}, \bx_1)$. The numbers $m_{\ell}$ can be very different for other QoIs.

\begin{table}[ht!]
\centering
\begin{tabular}{|c|c|c|c|c|c|} \hline
  $\varepsilon$ & $m_0$ & $m_1$ & $m_2$ & $m_3$ & $m_4$ \\ \hline
0.2 & 28 &0 &0 &0 &0 \\ \hline
0.1 & 10 & 1 &0 &0 &0 \\ \hline
0.05 & 57 & 4 & 1 & 0& 0\\ \hline
0.025 & 342 & 22 & 4 & 1 &0\\ \hline
0.01  & 3278 & 205 & 35 & 6 & 1\\ \hline
\end{tabular}
\captionsetup{width=.75\textwidth}
\caption{The number of samples $m_{\ell}$ vs. $\varepsilon$, QoI $g=\sol_m(t_{15}, \bx_1)$.}
 \label{tab:ml}
\end{table}
\subsection{Comparison of MC and MLMC methods}
\begin{figure}[t!]
\begin{center}
\includegraphics[width=0.48\textwidth]{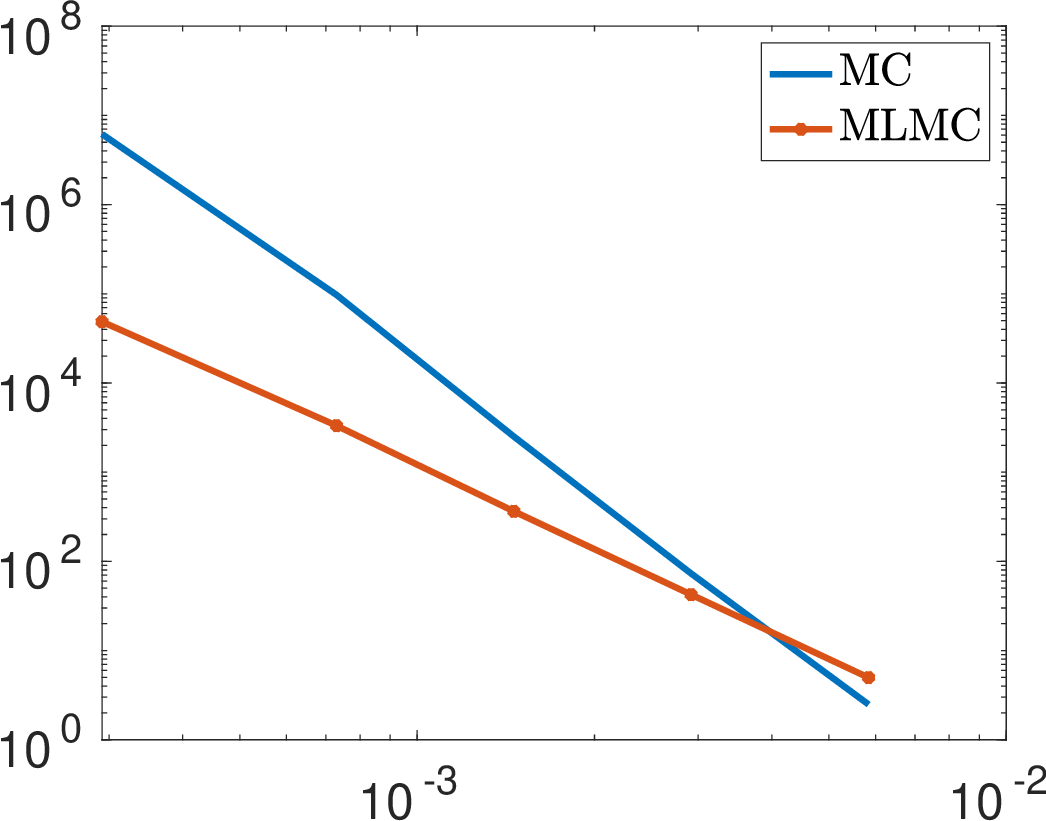}
\includegraphics[width=0.50\textwidth]{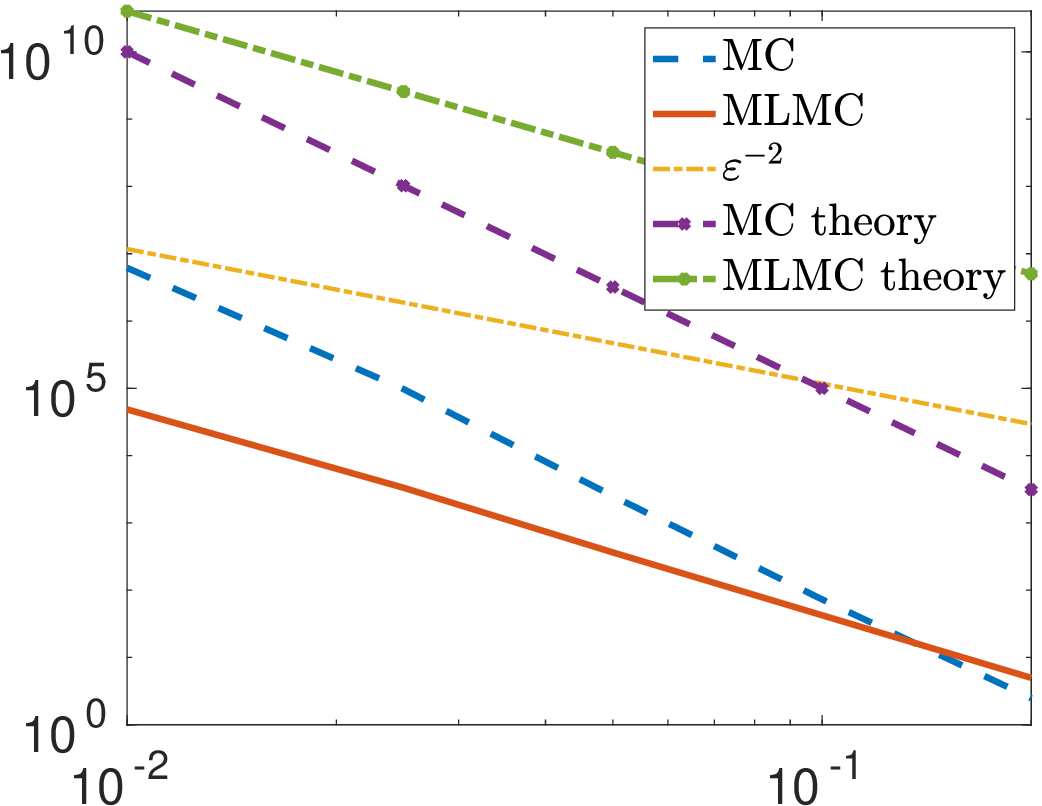}
\caption{(left) Computational cost comparison of MC and MLMC methods against $\varepsilon$ (horizontal axis) in log-log scale. (right) Same as on the left with additional graphs showing the slopes of the theoretical costs.
The cost (in seconds) is on the $y$ axis.}
\label{fig:MLMCvsMC}
\end{center}    
\end{figure}

In Figure~\ref{fig:MLMCvsMC}(left) we compare the estimated costs of the MC (denoted by a blue solid line) and MLMC (red solid line) methods. The cost of MLMC is calculated as in \eqref{eq:total_cost_MLMC}. We only give the estimated cost of MC because it is too time-consuming to compute. Such a comparison is typical for MLMC papers \cite{CMLMC,haji2015multi,ErikOptGeom15}. In addition, in Figure~\ref{fig:MLMCvsMC}(right) we plot the slopes of the $\varepsilon^2$ graph (yellow solid line) and the theoretical MC and MLMC costs, purple dotted line and red dashed line, respectively. 
Note that these three graphs are not scaled to make them easier to read, only the slope is important.
The theoretical MC cost is $\mathcal{O}(\varepsilon^{-2-\frac{\hat{d}\gamma}{\alpha}})$, with $\hat{d}=2+1=3$, $\gamma=1$ (linear cost of the geometric multigrid method), $\alpha=1.08$. Thus, the theoretical MC cost is $\mathcal{O}(\varepsilon^{-5})$.
The cost of MLMC, according to the third case in \eqref{eq:mlmc_iso_work}, where $\beta<\tilde{d}\cdot \gamma$, is $\mathcal{O}(\varepsilon^{-(2+\frac{3\cdot 1 - 2}{1})})=\mathcal{O}(\varepsilon^{-3})$. This graph is also up to a constant factor.
From the left figure, we can see that the MLMC outperforms the MC method for almost all values of $\varepsilon$. The right figure also shows that the theoretical MC and MLMC lines (two top plots) have the same slope as the computed MC and MLMC lines (two bottom plots).

\newpage
\section{Discussion and Conclusion}
\label{sec:Conclusion}

In this study we examined a setting that mimics the Henry-like problem \cite{Simpson2003,Simpson04_Henry}, modeling seawater intrusion into a 2D coastal aquifer. The solution shows how pure water recharge (see Fig.~\ref{fig:Henry2d-scheme}(left)) from the ``land side" mitigates the salinization of the aquifer due to the influx of saline water from the ``sea side," thereby achieving an equilibrium in salt concentration. Building on \cite{GRILLO2010}, we incorporated a fracture on the seaside that significantly enhances the permeability of the porous medium.

Uncertainties were introduced in the fracture width, the porosity of the bulk medium, its permeability and the pure water recharge from the land side. Porosity and permeability were modelled as random fields, recharge as a random periodic intensity and fracture width as a random variable. In total, three random variables were involved. We calculated the mean and variance of the salt mass fraction, which is subject to uncertainty.

The central question of this study was to evaluate the efficiency of the Multi-Level Monte Carlo (MLMC) method in computing statistics for different quantities of interest (QoIs). Our results indicate that the efficiency of MLMC depends on the QoI $g_{\ell}$, on the convergence rates of $\mathbb{E}[g_{\ell} - g_{\ell-1}]$ and $\Var{g_{\ell} - g_{\ell-1}}$ (see Table~\ref{tab:estimated_rates} and Fig.~\ref{fig:MLMC_MC_mean}), which should remain consistent across different levels $\ell$. In our numerical tests, which we do not publish here, we found that for some QoIs the convergence rates change, i.e. do not remain consistent across levels $\ell = 1, \ldots, 5$. One such example is the integral $I_4(\sol_m)$, as given in \eqref{eq:integral_box}. In particular, the rate $\alpha$ for levels $\ell = 1, \ldots, 4$ was equal to $-1.5$, but for $\ell = 5$ it was a completely different number. A possible explanation is that the 5th level is not needed. This observation needs further investigation.

Another challenge was the time dependence. The number of levels $L$ and the number of samples $m_{\ell}$ depend on the decay of the mean and variance (see \eqref{eq:num_levels}-\eqref{eq:scaled_m_l_S} and Figures~\ref{fig:QOIs_for_point1} bottom left and bottom right), which in turn depend on time $t$. Initially, when $t$ was small, the variability was low, increased during the mixing process of salt and fresh water, and decreased when equilibrium was reached. So, it is impossible to provide $L$ and $m_{\ell}$, which are optimal for every time step.

The number of samples required at each level was estimated by minimising the objective function (see \eqref{eq:goal_function}) used in the MLMC algorithm. A different objective function has been used in some other MLMC work (see \cite{CMLMC}). The choice of the goal function can affect the MLMC efficiency.

To achieve the efficiency of the MLMC approach demonstrated in this study, it is crucial that the complexity of solving each random realization scales with the number of grid vertices at each level. This was achieved using the geometric multigrid method implemented in the ug4 software toolkit \cite{ug4_ref1_2013, ug4_ref2_2013}, and our numerical tests validated the expected efficiency. The random realizations were computed concurrently on different nodes of the Shaheen II supercomputer at the King Abdullah University of Science and Technology.

Our results confirm that the MLMC method is applicable to the problem at hand, computing the mean up to 100 times faster than standard Monte Carlo methods (see Fig.~\ref{fig:MLMCvsMC}, left). Sampling across different mesh levels effectively reduces the overall computational cost.\\
\textbf{Constraints.}
Time dependence poses challenges, as the optimal number of samples varies with time $t$ and spatial location $\mathbf{x}$, potentially requiring fewer samples at some points and more at others.
Each new QoI necessitates re-evaluation of parameters $\alpha$ and $\beta$, resulting in different MLMC performance characteristics for each QoI.\\
\textbf{Future Work.} Further research could focus on more realistic modeling of fractures, porosity, and permeability by incorporating additional random variables.
One could allow more fractures, their uncertain orientations.
One can further research how geological parameters, including the fracture, influence flow and transport.

\backmatter

\bmhead{Acknowledgments}
For computing time, this research used Shaheen II, which is managed by the Supercomputing Core Laboratory at the King Abdullah University of Science and Technology (KAUST) in Thuwal, Saudi Arabia. We thank the KAUST HPC support team for their assistance with Shaheen II. This work was supported by the Alexander von Humboldt foundation. 
%






\bibliography{MoCa}

\end{document}